\documentclass[12pt]{article}
\usepackage{amsmath}
\usepackage{amsthm}
\usepackage{amssymb}
\usepackage{stmaryrd}
\usepackage{graphicx,psfrag,epsf}
\usepackage{enumerate}
\usepackage{natbib}
\usepackage{hyperref}
\usepackage{array}
\usepackage{multirow}

\newcommand{\blind}{1}

\newtheorem{prop}{Proposition}
\newtheorem{defn}{Definition}
\newtheorem{assumption}{Assumption}

\begin{document}

\def\spacingset#1{\renewcommand{\baselinestretch}%
{#1}\small\normalsize} \spacingset{1}


\if1\blind
{
  \title{\bf  Optimal control   for parameter estimation in partially observed  hypoelliptic stochastic differential equations}
  \author{Quentin Clairon  
\hspace{.2cm}\\
  University of Bordeaux, Inria Bordeaux Sud-Ouest, \\
Inserm, Bordeaux Population Health Research Center, SISTM Team, UMR1219, \\
F-33000 Bordeaux, France\\
Vaccine Research Institute, F-94000 Créteil, France\\
    and \\
    Adeline Samson  \thanks{
    The author gratefully acknowledges the LabEx PERSYVAL-Lab
(ANR-11-LABX-0025-01) and  MIAI@Grenoble Alpes, (ANR-19-P3IA-0003).}\\
    Univ. Grenoble Alpes, CNRS, Inria, Grenoble   Institute of Engineering,  LJK,\\
  38000 Grenoble, France
}
  \maketitle
} \fi

\if0\blind
{
  \bigskip
  \bigskip
  \bigskip
  \begin{center}
    {\LARGE\bf  Optimal control  for parameter estimation in partially observed  hypoelliptic stochastic differential equations}
\end{center}
  \medskip
} \fi

\bigskip
\begin{abstract}
We deal with the problem of parameter estimation in stochastic differential equations (SDEs) in a partially observed framework. We aim to design a method working for both elliptic and hypoelliptic SDEs, the latters being characterized by  degenerate diffusion coefficients. This feature often causes the failure of contrast estimators based on Euler Maruyama discretization scheme and dramatically impairs classic stochastic filtering methods used to reconstruct the unobserved states. All of theses issues make the estimation problem in hypoelliptic SDEs difficult to solve. To overcome this, we construct a well-defined cost function no matter the elliptic nature of the SDEs. We also bypass the filtering step by considering a control theory perspective. The  unobserved states are estimated by solving deterministic optimal control problems using  numerical methods which do not need strong assumptions on the diffusion coefficient conditioning. Numerical simulations made on different partially observed hypoelliptic SDEs reveal our method produces accurate  estimate while dramatically reducing the computational price comparing to other methods. 
\end{abstract}

\noindent%
{\it Keywords:}  Stochastic differential equations, Parameter estimation, Hypoellipticity, Optimal control theory
\vfill

\newpage
\spacingset{1.45} 

\section{Introduction\label{sec:introduction}}

Stochastic differential equations (SDEs) are widely used in several fields to model stochastic temporal dynamics. We focus on hypoelliptic SDE, that is when the diffusion matrix is not of full rank but its solution has a smooth density. More precisely, we consider a process   $Z_{t}=(V_t, U_t)\in\mathbb{R}^{d}$,  of dimension $d$   on a time interval $[0,T]$, defined as the solution of the following hypoelliptic SDE:
\begin{equation}
\begin{array}{l}
dV_{t}=g_{\theta}(V_{t},U_{t},t)dt\\
dU_{t}=h_{\theta}(V_{t},U_{t},t)dt+B_{\sigma}(Z_{t},t)dW_{t}\\
Z(0)=Z_{0}
\end{array}\label{eq:sde_original_model}
\end{equation}
where $V_{t}$
is the $d_{V}$-dimensional component corresponding to the smooth
part of the process, that is no directly affected by stochastic perturbations,
and $U_{t}$ is the  $d_{U}$-dimensional rough part, 
$W_{t}$ is a $d_{U}$ dimensional Brownian motion acting on
the system through the $d_{U}$ squared diffusion matrix
$B_{\sigma}(Z_{t},t)$. In some applications, only the smooth component can be measured.   

Let us give some examples of such models. Stochastic damping Hamiltonian systems or Langevin equations have been developed to describe  a particle moving in an environment defined by a potential \cite{Wu2001}.  The system is usually two dimensional, the first equation describes the position  of the particle and the other its velocity. The noise is degenerate because it acts only on the velocity and not on the position.  Only the position is directly measured. We  can cite applications of these models in molecular dynamics \cite{LeimkuhlerMatthewsBook2015}, paleoclimate research \cite{Ditlevsen2002},   neural field models \cite{CoombesByrne2017,DitlevsenLocherbach2017}.  Another application is neuronal models of membrane potential dynamics. Examples are the hypoelliptic FitzHugh-Nagumo (FHN) model
\cite{Gerstner2002, DeVille2005,LeonSamson2017}, the hypoelliptic Hodgkin-Huxley model 
\cite{GoldwynSheaBrown2011,TuckwellDitlevsen2016}, or synaptic-conductance based models with
stochastic channel dynamics \cite{Paninski2012,DitlevsenGreenwood2013}. Only the first equation, corresponding to the membrane potential, can be measured by intra-cellular recordings.   It is therefore
important to develop   estimation methods for this class of
models.

In the rest of paper, we assume the drift
functions $g_{\theta}$ and $h_{\theta}$ depend on a $d_{\theta}$
dimensional parameter $\theta$ and the matrix $B_{\sigma}$ on 
the $d_{\sigma}$ dimensional parameter $\sigma$ known as volatility. 
Partial observations of $(Z_t)$ are defined as:
\[
Y_{t}=CZ_{t}
\]
where $C$ is the $d_{C}\times d$   observation matrix. We
assume $(Y_t)$ is discretely observed on the interval $\left[0,\,T\right]$
at times $0=t_{0}<\ldots<t_{n}=T$ and without measurement error.
Our aim is to estimate the unknown parameter $\psi=(\theta,\sigma)$
of model (\ref{eq:sde_original_model}) using the  discrete and partial
observations $(Y_{0},\ldots,Y_{n})$ and possibly without knowing
the initial condition $Z_{0}.$

 Parametric estimation of hypoelliptic models faces several difficulties. When the first equation reduces to $g_\theta(V_t,U_t)=U_t$, different contrasts estimators haven been proposed \cite{genon-catalot2000,Gloter2000, Ditlevsen2004,Gloter2006,Samson2012, Lu2016, Leon2019}.  This specific case allows to deal with partial observations, the second coordinate $U_t$ being replaced by increments of $V_t$. When the drift $g_\theta$ is more complex, we deal with several entangled issues. 
 First,  solutions of the SDE are generally non explicit.  
 Numerical schemes are used to approximate the discretized process and derive an estimation criteria. However, compared to elliptic system, the degeneracy of the noise complicates the statistical estimation. Thus estimation cost function   that would   be derived directly from the Euler-Maruyama fails because  the covariance matrix $\varGamma_{\sigma}^{T}\varGamma_{\sigma}$ where $\varGamma_{\sigma}(Z_{t},t)=\left(0_{d_{V},d_{U}}, 
B_{\sigma}(Z_{t},t)\right)^t$, is not   full rank. \cite{Gloter2021} propose a local Gaussian approximation of the transition density in the case of totally observed coordinates. Several papers propose to use a likelihood based on higher order schemes of approximation,  the 1.5 order scheme \cite{DitlevsenSamson2019}, the Local Linearization scheme \cite{Melnykova2020}. 
In the case of partial observations, the unobserved coordinates $U_t$ has to be filtered or imputed, either in a Bayesian settings with a Markov Chain Monte Carlo algorithm   \cite{Pokern2009,Graham2019} or in a frequentist settings with a particle filter or sequential Monte Carlo coupled with a Stochastic Approximation
Expectation Maximization (SAEM) algorithm  \cite{Ditlevsen2014, DitlevsenSamson2019}. A recent simulation filter has been proposed by \cite{Bierkens2020} for hypoelliptic diffusions. 
Still the theoretical assessment and
practical efficiency of such numerical methods critically rely on
$\varGamma_{\sigma}$ being well conditioned. These methods are also time consuming.

In this work, we propose
an estimation method  based on optimal control theory which applies to nonlinear hypoelliptic systems and partial observations. The idea of using optimal control theory is to treat the coordinate filtering problem  as a deterministic tracking problem. Then, numerical devices coming from control theory allow us to estimate the \textit{optimal control} bringing the observed system states  the closest possible to the observations.  This was proposed for ODE system  
\cite{BrunelClairon_Pontryagin2017,clairon2019tracking,clairon2021regularization, Iolov2017}.  
For SDE, the only reference is \cite{clairon2020optimal} which uses linear-quadratic theory, a particular case of the Pontryagin maximum principle, to solve the tracking problem. But this limits the application of the method to linear systems. Also, this method separatly estimates $\sigma$ and $\theta$ via a nested procedure because of identifiability issues. Volatility $\sigma$ was estimated via the minimization of an outer criteria demanding repeated estimation of $\theta$ based on the optimization of an inner criteria. This dramatically increases the computation time of the method when $\sigma$ is unknown. Additionaly, the criteria used for $\sigma$ inference was properly defined only for linear SDEs. 

In this work, we construct an estimation procedure which simultaneously estimate $(\theta, \sigma)$ by relying on a contrast function bypassing identifiability issues for $\sigma$ estimation encountered in  \cite{clairon2020optimal}. This contrast function is well defined for linear and nonlinear SDEs, both elliptics and hypoelliptics. Also, the control problem solutions required by our approach are easily obtained  for non-linear SDEs  via an adaptation of the method proposed by \cite{CimenBanks2004}.

Our statistical criteria exploits the hypoelliptic nature of (\ref{eq:sde_original_model})
which assumes  "enough" interactions between $V_{t}$ and $U_{t}$
such that  the Wiener process   $W_{t}$  
influences all the coordinates of the SDE, potentially with some delay. More precisely, we consider a lagged
contrast function linking $W_{t}$ and the observed part of the system
state $CZ_{t+m_{B}}$ after a delay $m_{B}>0$ such that the dependence
$CZ_{t+m_{B}}$with respect to $W_{t}$ is regular enough to ensure
 a well-conditioned variance for $CZ_{t+m_{B}}$. This   requires to know explicitly 
 the dependence of $CZ_{t+m_{B}}$   to
$W_{t}$. For this we rely on the  so-called \textbf{pseudo-linear} representation of the discretized  SDE (\ref{eq:sde_original_model}) \cite{Cimen2008}. The scheme is expressed at each time $t$ as a function of the $m_B$ previous lagged states. The discretized process is no more Markovian, this is a direct consequence of hypoellipticity. From this lagged formulation of the discretized process, we deduce a  lagged contrast which is the pseudo-likelihood of the representation. This statistical contrast depends on  the  whole state variable, even the unobserved ones. We thus need to filter or predict the unobserved coordinate $U_t$. 
We propose a predictor which is a   balance between
fidelity to the observations and to the solution of SDE (\ref{eq:sde_original_model}).
This   predictor is defined via 
an optimal control problem that can be solved easily without assumption
on $\varGamma_{\sigma}$ rank. The optimal control problem is solved
with the linear quadratic theory for linear SDEs and for nonlinear
models with an adaptation of the procedure developed by \cite{Cimen2004}
  based on the pseudo-linear representation. 
  Our criteria is somewhat similar to the Generalized profiling introduced
in \cite{Ramsay2007} for Ordinary differential equations but here
dedicated to SDEs with potentially degenerate diffusions. 
We also consider a criteria making a tradeoff between data and model
fidelity to regularize the inverse problem of parameter estimation.
This method provides estimators for the parameters and for the unobserved coordinate $U_t$. The simulation study illustrates the good performances of our method.

Section \ref{sec:model} studies the hypoellipticity of the continuous process, provides a discrete approximation scheme and the corresponding concept of hypoellipticity. 
In section \ref{sec:lagged_contrast},  we  derive the expression
of our statistical criteria which defines our estimator. The construction
of the   state predictor and  the optimal control numerical methods   are presented in section \ref{sec:state_predictor_estimation}.
Then, in section \ref{sec:Simulation} we conduct a simulation study
on   three hypoelliptic models to evaluate the practical accuracy
of our method as well as its computational efficiency.

\section{Model and its discrete approximation \label{sec:model}}
This Section describes the properties of hypoelliptic systems and the connexity property that allows to propagate the noise to the smooth coordinates. Then we introduce the approximate Euler-Maruyama scheme and the lagged pseudo-linear formulation, the key of our estimation method. 
 \subsection{Models and its assumptions}

Let us denote: 
\[
f_{\theta}(Z_{t},t)=\left(\begin{array}{c}
g_{\theta}(V_{t},U_{t},t)\\
h_{\theta}(V_{t},U_{t},t)
\end{array}\right),\,\varGamma_{\sigma}(Z_{t},t)=\left(\begin{array}{c}
0_{d_{V},d_{U}}\\
B_{\sigma}(Z_{t},t)
\end{array}\right)
\]
 the   drift and   diffusion coefficient of the diffusion process $Z_t$ respectively, then  
$$dZ_{t}=f_{\theta}(Z_{t},t)dt+\varGamma_{\sigma}(Z_{t},t)dW_{t}.$$
We assume $B_{\sigma}$ is   full rank and we consider two classes of models called elliptic and hypoelliptic. 

\vspace{0.5em}
  \textbf{Elliptic models} are defined by $d_{V}=0$. In that case,  the Brownian process directly acts on the whole
system, the diffusion coefficient is not degenerated. Girsanov formula applies and the process has a continuous smooth density.

\vspace{0.5em}
 \textbf{Hypoelliptic models} are defined with $d_{V}>0$  and verify the H\"{o}rmander condition detailed below. 
Let us first  defined the  Lie brackets. 
For $f$ a function with values in $\mathbb{R}^{d}$, $f_{l}$
stands for its $l$-th component.  For a matrix $\varGamma$,  $\varGamma_{j}$ denotes its $j-$th column. 

\begin{defn}
The Lie bracket of two functions $f,g:\mathbb{R}^{d}\rightarrow \mathbb{R}$ is defined as 
$$\left[f,\,g\right]_{l}=\frac{\partial g_{l}(z)}{\partial z}^{T}f(z)-\frac{\partial f_{l}(z)}{\partial z}^{T}g(z).$$
\end{defn}

Hypoellipticity is then defined as follows:

\begin{defn} 
(H\"{o}rmander condition)
Let us consider  the set   $\mathcal{L}$     defined iteratively as: 
\begin{itemize}
\item Initialization step:  $\mathcal{L}$ is composed of vectors $L^0=\varGamma_{j}$, for  $j=1,\ldots d_{U}$.
\item Generalization step at iteration $k$,  if $L^{k-1}\in\mathcal{L}$ then vectors $L^k$ defined by $L^k=\left[f(Z_{t},t),\,L^{k-1}\right]$ and by   $L^k=\left[\varGamma_{j}(Z_{t},t),\,L^{k-1}\right]$
for $j=1,\ldots d_{U}$ belong  to $\mathcal{L}$.
\end{itemize}
If at some iteration, $\mathcal{L}$ spans
$\mathbb{R}^{d}$,  the weak H\"{o}rmander condition is fulfilled and the  system is said hypoelliptic.
\end{defn}

\subsection{Hypoellipticity properties}

We derive a necessary
condition on the drift $g_{\theta}$ respected by hypoelliptic systems
(the proof is given in Appendix 1 in Supplementary Materials):
\begin{prop} (Connexity property)
\label{prop:connexity_prop}If the SDE (\ref{eq:sde_original_model})
is hypoelliptic then for each $l\in\left\llbracket 1,d_{V}\right\rrbracket $,
it exists $j\in\left\llbracket 1,d_{U}\right\rrbracket $ and a sequence
$\left(q_{1},\ldots q_{m_{l}}\right)\in\left\llbracket 1,\,d_{V}\right\rrbracket ^{m_{l}}$
with $m_{l}\leq d_{V}$ such that 
\begin{equation}
\frac{\partial g_{\theta,l}(Z_{t},t)}{\partial V_{q_{1}}}\frac{\partial g_{\theta,q_{1}}(Z_{t},t)}{\partial V_{q_{2}}}\cdots\frac{\partial g_{\theta, q_{m_{l}-1}}(Z_{t},t)}{\partial V_{q_{m_{l}}}}\frac{\partial g_{\theta,q_{m_{l}}}(Z_{t},t)}{\partial U_{j}}\neq0.\label{eq:connexity_prop}
\end{equation}
\end{prop}

This can be seen as a connexity property between the set of smooth
variables and rough ones. It means that each variable $V_{l,t}$ is
influenced by at least one component of $U_{t}$ at least undirectly
through a path of auxiliary smooth variables ($V_{q_{1}},\ldots, V_{q_{m_{l}}}$).
Of course, from a given sequence $(q_1, \ldots, q_{m_l})$, it is possible to construct bigger
ones. To remove ambiguity $m_{l}$   denotes the length of the
shortest  possible path   such that the connexity
property (\ref{eq:connexity_prop}) holds.

 Although the diffusion coefficient $\varGamma_{\sigma}(Z_{t},t)$ is rank
deficient, the solution of (\ref{eq:sde_original_model})   has
a smooth density with respect to the Lebesgue measure.

If the drift $f(Z_t)$ satisfies a dissipativity condition
$$ <f(z),z>\leq \alpha -\delta \|z\|^2,\forall z \in \mathbb{R}^d $$
where $\alpha, \delta>0$, the function  $L(z) = 1+ \|z\|^2$ is a Lyapounov function. It ensures the geometric ergodicity of the solution of the SDE. It means that $(Z(t))_{t \in [0,T]}$ converges exponentially fast to a unique invariant distribution. 

However, usually, the solution of  SDE (\ref{eq:sde_original_model}) is not explicit and a numerical scheme is needed to approximate the solution.

\subsection{Euler-Maruyama discretization}

 Let us introduce the Euler-Maruyama discretized process ($\tilde Z_{i}$, for $i=1, \ldots, n$) of SDE (\ref{eq:sde_original_model}). To simplify the notations, we assume that the time points are equidistant, that is $\triangle_{i}=t_{i+1}-t_{i}= \triangle, \forall i=1, \ldots, n$. Then the discretized approximate process is defined by: 
\begin{equation}
\begin{array}{l}
\tilde Z_{i+1}=\tilde Z_{i}+\triangle f_{\theta}(\tilde Z_{i},t_{i})+\sqrt{\triangle}\varGamma_{\sigma}(\tilde Z_{i},t_{i})u_{i}\\
\tilde Z_{0}=Z_{0}
\end{array}\label{eq:discretized_original_model}
\end{equation}
where    $u_{i}\sim N(0,I_{m})$ is the normalized increment of the Brownian
motion at $t_{i}$, i.e. $u_{i}=\frac{1}{\sqrt{\triangle }}\left(W_{t_{i+1}}-W_{t_{i}}\right)$ and 
   $u=\left(u_{0},\ldots,u_{n-1}\right)$. 
   
   The connexity property can be illustrated on this Euler-Maruyama discretisation.
Proposition \ref{prop_connexity_numerical_scheme} proves that  $m_{l}+1$ iterations of (\ref{eq:discretized_original_model}) are required to have $\tilde V_{l,i+m_{l}+1}$
affected by some entry $\tilde U_{j,i}$ of the rough component.  
\begin{prop}
\label{prop_connexity_numerical_scheme} Let us assume that for each
$l\in\left\llbracket 1,d_{_{V}}\right\rrbracket $ it exists a unique
minimal length sequence $\left(q_{1},\ldots q_{m_{l}}\right)$ such
that the connexity property   holds. Then
the Euler-Maruyama approximation ($\tilde Z_{i}$, for $i=1, \ldots, n$) is such that for each $l\in\left\llbracket 1,d_{_{V}}\right\rrbracket$ it exists $j\in\left\llbracket 1,d_{U}\right\rrbracket $ such that:
 $$\frac{\partial \tilde V_{l,i+k}}{\partial \tilde U_{j,i}}=0, \quad \mbox{for } k=1, \ldots, m_l \quad \mbox{and}\quad \frac{\partial \tilde V_{l,i+m_{l}+1}}{\partial \tilde U_{j,i}}\neq0$$
\end{prop}

So, in model (\ref{eq:discretized_original_model}), each coordinate 
of $\tilde V_{i+m_{B}+1}$ with $m_{B}:=\max_{l}{m_{l}}$ is influenced
by $\tilde U_{i}$. It also means that $m_{B}+1$ is the minimal delay  to propagate the noise $u_i$ from the rough components to all the smooth ones in the discretized model (\ref{eq:discretized_original_model}).

Let us introduce a new observability assumption that reinforces  hypoellipticity for a discrete process. It requires that the observed variables give enough information about the diffusion process, or equivalently that the $Y_t$ are influenced by all the rough components.
\begin{assumption}
(H1): The $d_{o}\times d_{U}$ matrix
\[
C\left(\prod_{l=1}^{m_{B}}\mathbf{A}_{\theta,i+l}\right)\varGamma_{\sigma}(z_i,t_{i})
\]
is   full rank for all $1\leq i\leq n-m_{B}$ and for all possible values
$\left\{ z_i\right\} _{0\leq l\leq m_{B}}\in\mathbb{R}^{d\times m_{B}}.$
\end{assumption}

Let us detail this assumption for elliptic and some hypoelliptic models. 

\begin{itemize}
\item Elliptic models verify   $m_{B}=0$.  Assumption (H1) imposes
that $C\varGamma_{\sigma}(\tilde Z_{i},t_{i})$ is   full rank. It ensures a smooth transition density for the discrete process ($Y_i$). This is similar to the continuous process. 

\item We focus on hypoelliptic SDE with $\frac{\partial g_{\theta}}{\partial V}(V_{t},U_{t},t)$
  full rank and $Y_t=V_t$.  This corresponds to the case $m_{B}=1$, that is all the components of
$V_{t}$ have   a direct interaction with at least one component
of $U_{t}$ and   stochasticity only acts on the unobserved variables.
  This corresponds to the 1-step hypoellipticity defined by \cite{buckwar2021splitting}. 

\end{itemize}

\subsection{Lagged pseudo-linear formulation of the discretized process}
Proposition \ref{prop_connexity_numerical_scheme} proves 
the existence of a minimal delay $m_{B}$. For estimation purpose,
we   need the explicit link between $\tilde Z_{i+m_{B}+1}$
and $u_{i}$. For this, we rely on what we call later the \textit{lagged  pseudo-linear formulation} of
(\ref{eq:discretized_original_model}). 

Let us  introduce the matrix
$A_{\theta}(\tilde Z_{i},t_{i})$ and the vector $r_{\theta}(t_{i})$
such that $A_{\theta}(\tilde Z_{i},t_{i})\tilde Z_{i}+ r_{\theta}(t_{i})=f_{\theta}(\tilde Z_{i},t_{i})$. This allows  to reformulate (\ref{eq:discretized_original_model}) with a pseudo-linear expression:
\begin{equation}
\left\{ \begin{array}{lll}
\tilde Z_{i+m_{B}+1} & = & \mathbf{A}_{\theta,i+m_{B}}\, \tilde Z_{i+m_{B}}+\triangle r_{\theta}(t_{i+m_{B}})\\
 & + & \sqrt{\triangle }\, \varGamma_{\sigma}(\tilde Z_{i+m_{B}},t_{i+m_{B}})\, u_{i+m_{B}}\\
\tilde Z_{0} & = & Z_{0}
\end{array}\right.\label{eq:sde_optimal_control_notation}
\end{equation}
with $\mathbf{A}_{\theta,i}:=I_{d}+\triangle A_{\theta}(\tilde Z_{i},t_{i})$. 
We  can replace $\tilde Z_{i+m_{B}}$ by its state-space expression (\ref{eq:discretized_original_model}) and iterate to obtain the lagged pseudo-linear formulation.

\begin{prop} 
The lagged pseudo-linear formulation of the discretized process
(\ref{eq:discretized_original_model}) is defined as follows. 
For $i=0,\ldots,n-m_{B}-1$:
\begin{equation}
\begin{array}{lll}
\tilde Z_{i+m_{B}+1} & = & \prod_{r=0}^{m_{B}}\mathbf{A}_{\theta,i+r}\tilde Z_{i}\\
 & + & \sum_{r=0}^{m_{B}}\left(\prod_{l=r+1}^{m_{B}}\mathbf{A}_{\theta,i+l}\right)\triangle r_{\theta}(t_{i+r})\\
 & + & \sum_{r=0}^{m_{B}}\left(\prod_{l=r+1}^{m_{B}}\mathbf{A}_{\theta,i+l}\right)\sqrt{\triangle }\varGamma_{\sigma}(\tilde Z_{i+r},t_{i+r})u_{i+r}.
\end{array}\label{eq:Lagged_recursive_formulation}
\end{equation}
In
the last expression the multiplication symbol $\prod_{l=0}^{m_{B}-1}$
 denotes the sequential left multiplication i.e.
$\prod_{j=a}^{b}\mathbf{A}_{\theta,j}=\mathbf{A}_{\theta,b}\times\mathbf{A}_{\theta,b-1}\times\ldots\times\mathbf{A}_{\theta,a}$
with the convention $\prod_{j=b+1}^{b}\mathbf{A}_{\theta,j}=I_{d}.$
\end{prop}

The pseudo-linear formulation gives the explicit link between $\tilde Z_{i+m_{B}+1}$ and $u_{i}$. 
The Markovian property is no more fulfilled for the discretised process: the smooth density is obtained only if we express $\tilde Z_{i+m_{B}+1}$ as a function of the $m_B$ previous states. This is a direct consequence of hypoellipticity.

\section{Cost function construction \label{sec:lagged_contrast}}

Our aim is to estimate the unknown parameter $\psi=(\theta,\sigma)$
of model (\ref{eq:sde_original_model}) using the discrete and partial
observations $(Y_{0},\ldots,Y_{n})$ and possibly without knowing
the initial condition $Z_{0}.$

We first assume that the process $Z$ is observed and construct an estimation contrast based on the lagged pseudo-linear representation. Then, we propose to filter $Z$ by optimal control. 

\subsection{Statistical criteria assuming $Z$ observed}
Let us consider in this Section that the variables $Z$ are observed.  We can use formula (\ref{eq:Lagged_recursive_formulation})   to build a statistical criteria. 
Let us denote $t_{i+m_{B}}^{i}=\left(t_{i},\ldots,t_{i+m_{B}}\right)$,  $\tilde Z_{i:i+m_{B}}=\left(\tilde Z_{i}\ldots,\tilde Z_{i+m_{B}}\right)$ and define 
\[
\begin{array}{l}
X_{i:i+m_{B}+1}=Y_{i+m_{B}+1}-C\left(\prod_{r=0}^{m_{B}}\mathbf{A}_{\theta,i+r}\right)\tilde Z_{i}-\sum_{r=0}^{m_{B}}F_{i+r:i+m_{B}}\\
F_{i+r:i+m_{B}}=\triangle C\left(\prod_{l=r+1}^{m_{B}}\mathbf{A}_{\theta,i+l}\right)r_{\theta}(t_{i+r})\\
G_{i+r:i+m_{B}}=\sqrt{\triangle }C\left(\prod_{l=r+1}^{m_{B}}\mathbf{A}_{\theta,i+l}\right)\varGamma_{\sigma}(Z_{i+r}^{d},t_{i+r}).
\end{array}
\]
Then we have
\begin{equation}
\begin{array}{l}
X_{i:i+m_{B}}=\sum_{r=0}^{m_{B}}G_{i+r:i+m_{B}}u_{i+r}.\end{array}\label{eq:value_Xi}
\end{equation}
Under assumption (H1), the matrix 
 $$\varSigma_{i:i+m_{B}}=\sum_{r=0}^{m_{B}}G_{i+r:i+m_{B}}G_{i+r:i+m_{B}}^{T}$$
 is nonsingular. Then, conditionnally on $\tilde Z_{i:i+m_{B}}$, we have 
\[
X_{i:i+m_{B}+1}\mid \tilde Z_{i:i+m_{B}}\sim N(0_{d_{o},1},\,\varSigma_{i:i+m_{B}})
\]

and we can define the following cost function based on the log likelihood:
\begin{equation}
H^{m_{B}}(\psi\mid Y, Z) = -log \left(\mathbb{P}_{\psi}\left\{ X_{i:i+m_{B}+1}\mid \tilde Z_{i:i+m_{B}}\right\} _{i\in\left\llbracket 1,\,n- m_B - 1\right\rrbracket }\right)
\label{eq:H_mB}
\end{equation} and the associated estimator: 
\begin{equation}
\begin{array}{lll}
\widehat{\psi}&=&\arg\min_{\psi}H^{m_{B}}(\psi\mid Y, Z)\label{eq:H_mB}.
\end{array}
\end{equation}
From equation (\ref{eq:value_Xi}), we see that $X_{i:i+m_{B}+1}$ depends at most of $\left(u_{i},\ldots,\,u_{i+m_{B}}\right)$. Thus the sequence  $(X_{i:i+m_B+1})$ for varying $i$ can be composed of dependant terms which makes $H^{m_{B}}$ not easy to calculate in this general setting.  We detail in the next subsection some cases where $H^{m_{B}}$ has an explicit form. 

\subsection{Simplification of $H^{m_{B}}$  for some models}
Note first that when the sequence of $(X_{i:i+m_B+1})$ is composed of independent terms, the expression of $H^{m_{B}}$ simplifies to 
\begin{equation}
 H^{m_{B}}(\psi\mid Y, Z) =   \sum_{i=1}^{n-m_{B}-1}\left(X_{i:i+m_{B}+1}^{T}\varSigma_{i:i+m_{B}}^{-1}X_{i:i+m_{B}+1}+\log\left(\left|\varSigma_{i:i+m_{B}}\right|\right)\right) .
\label{eq:H_mB2}
 \end{equation}
This is the case for elliptic and 1-step hypoelliptic systems. Let us give more
detailed expressions for $\varSigma_{i:i+m_{B}}$ for these two  cases.

\subsubsection{Elliptic system}
Here, we have $m_B=0$: 
\[
\begin{array}{l}
X_{i:i+1}=Y_{i+1}-C\left(\mathbf{A}_{\theta,i}\tilde Z_{i}+\triangle r_{\theta}(t_{i})\right)\end{array}
\]
and $X_{i:i+1}\mid \tilde Z_{i}\sim N(0_{d_{o},1},\,\varSigma_{i:i})\textrm{ with }\varSigma_{i:i}=\triangle C\varGamma_{\sigma}(\tilde Z_{i},t_{i})\varGamma_{\sigma}(\tilde Z_{i},t_{i})^{T}C^{T}.$
Thus, we fall back to the classic contrast estimator of elliptic SDE.

\subsubsection{1-step hypoelliptic system \label{subsec:1_step_hypoelliptic}}
In the case where $m_B=1$ and $V_{t}=CZ_{t}$, we end up with $C\varGamma_{\sigma}(\tilde Z_{i},t_{i})=0$. So condition
(H1) is respected if $CA_{\theta}(\tilde Z_{i+1},t_{i+1})\varGamma_{\sigma}(\tilde Z_{i},t_{i})$
is of  full rank. Then $X_{i:i+m_{B}+1}$ is given by:
\[
\begin{array}{l}
X_{i:i+2}=Y_{i+2}-C\mathbf{A}_{\theta,i+1}\mathbf{A}_{\theta,i}\end{array}\tilde Z_{i}-C\left(\triangle \mathbf{A}_{\theta,i+1}r_{\theta}(t_{i})+\triangle r_{\theta}(t_{i+1})\right)
\]
and follows the Gaussian law $X_{i:i+2}\mid \tilde  Z_{i:i+1}\sim N(0_{d_{o},1},\,\varSigma_{i:i+1})$
with:
\[
\varSigma_{i:i+1}=\triangle CA_{\theta}(\tilde Z_{i+1},t_{i+1})\varGamma_{\sigma}(\tilde Z_{i},t_{i})\varGamma_{\sigma}(\tilde Z_{i},t_{i})^{T}A_{\theta}(\tilde Z_{i+1},t_{i+1})C^{T}.
\]

\section{State predictor $\overline{Z}_{\theta,\sigma}$ \label{sec:state_predictor_estimation}}
Now, we propose a state predictor of the hidden process $Z$ using optimal control theory and then we plug it in the statistical criteria. 

\subsection{Formalization of the related optimal control problem}

In this section, we denote $(\tilde Z(Z_{0},u))$ the solution of (\ref{eq:sde_optimal_control_notation})
for a given   sequence  
$u=\left(u_{0},\ldots,u_{n-1}\right)$ and a given initial condition $Z_{0}$.
Similarly to   trajectory fitting estimators \cite{Kutoyants1991,Dietz2001},
a first   predictor  could be the solution $\overline{Z}$
of (\ref{eq:sde_optimal_control_notation}) which is the closest   to
the observations. This   solution corresponds to the sequence
$u_{Z_{0}}^{M}$ such that:
\begin{eqnarray}
u_{Z_{0}}^{M} & = & \arg\min_{u}\left\{ \sum_{i=0}^{n}\left\Vert C\tilde Z_{i}(Z_{0},u)-Y_{i}\right\Vert _{2}^{2}\right\} \label{eq:original_opt_control_definition}
\end{eqnarray}
when the initial condition $Z_{0}$ is fixed. However, this optimization
problem is ill-posed. The solution is not necessarily unique and would not always be continuous as a function of the parameters and observations.

We thus propose to introduce a penalization term into the optimization problem which will lead to a Tikhonov regularized version of it. This is known to remove sources of ill-posedness  \cite{Engl2009}, in particular it re-establishes  uniqueness and the continuity of the solution for linear SDEs. The 
 penalized  problem is the following:
\begin{eqnarray*}
\overline{u}_{Z_{0}} & := & \arg\min_{u}\left\{ \sum_{i=0}^{n}\left\Vert C\tilde Z_{i}(Z_{0},u)-Y_{i}\right\Vert _{2}^{2}-\frac{2}{w}\log P(u)\right\} 
\end{eqnarray*}
where $P(u)$ is the density of the increment $u$.
The term $-\frac{2}{w}\log P(u)$
penalizes the sequences which are unlikely to be a realization of the Brownian
motion and $w>0$ makes the balance between the  model
and data fidelity. We have
\[
\log P(u)=\log\prod_{i=0}^{n-1}P(u_{i})\varpropto\log\prod_{i=0}^{n-1}e^{-\frac{1}{2}u_{i}^{T}u_{i}}=-\sum_{i=0}^{n-1}\frac{1}{2}u_{i}^{T}u_{i}
\]
thus
$$
\overline{u}_{Z_{0}}   =  \arg\min_{u}\left\{ \sum_{j=1}^{n-1}\left(\left\Vert C\tilde Z_{i}(Z_{0},u)-Y_{i}\right\Vert _{2}^{2}+\frac{1}{w}u_{i}^{T}u_{i}\right)+\left\Vert C\tilde Z_{n}(Z_{0},u)-Y_{n}\right\Vert _{2}^{2}\right\}.
$$
The predictor corresponds to the solution $\overline{u}_{Z_{0}}$
of the following deterministic discrete optimal control problem:
\begin{equation}
\begin{array}{ll}
\textrm{Minimize in \ensuremath{u}:} & C_{w}(u|Y;Z_{0})=\sum_{j=1}^{n-1}\left(\left\Vert C\tilde Z_{i}(Z_{0},u)-Y_{i}\right\Vert _{2}^{2}+\frac{1}{w}u_{i}^{T}u_{i}\right)+\left\Vert C\tilde Z_{n}(Z_{0},u)-Y_{n}\right\Vert _{2}^{2}\\
\textrm{Subject to: } & \left\{ \begin{array}{l}
\tilde Z_{i+1}(Z_{0},u)=\mathbf{A}_{\theta,i}\tilde Z_{i}(Z_{0},u)+\Delta_{i}r_{\theta}(t_{i})+\sqrt{\triangle }\varGamma_{\sigma}(\tilde Z_{i}(Z_{0},u),t_{i})u_{i}\\
\tilde Z_{0}=Z_{0}.
\end{array}\right.
\end{array}\label{eq:OC_Cost_Regulator_Form}
\end{equation}
Thus, if we are able to solve problem (\ref{eq:OC_Cost_Regulator_Form})
to derive $\overline{u}_{Z_{0}}$, we can solve (\ref{eq:sde_optimal_control_notation})
to obtain the related state predictor $\overline{Z}_{Z_{0}}:=\tilde Z(Z_{0},\overline{u}_{Z_{0}})$. However,
our optimal control still depends on $Z_{0}$ which is potentially
unknown. If required,    its estimation is bypassed by profiling the cost
$C_{w}$ and   defining as initial condition estimator:
\[
\widehat{Z_{0}}=\arg\min_{Z_{0}}\left\{ \min_{u}C_{w}(u|Y;Z_{0})\right\}. 
\]
Let us now denote $\overline{u}:=\overline{u}_{\widehat{Z_{0}}}$ the corresponding
optimal control and $\overline{Z}:=\overline{Z}_{\widehat{Z_{0}}}$
the related state-space predictor. 

We propose   to solve (\ref{eq:OC_Cost_Regulator_Form})
with  numerical    control theory methods. For linear
models,   (\ref{eq:OC_Cost_Regulator_Form}) is a Linear-Quadratic
problem. The control theory then ensures:
\begin{itemize}
\item existence and uniqueness of the solutions $\overline{u}$ and $\overline{Z}$, which are continuous functions of   $\left(\theta,\sigma\right)$
and $Y$, under mild regularity hypothesis on $A_{\theta}$ and $\varGamma_{\sigma}$
(in particular, $\varGamma_{\sigma}$ is not required to be of full
rank),
\item computation of $\overline{u}$, $\overline{Z}$  and, if required, $\widehat{Z_{0}}$ by solving a finite
difference equation.
\end{itemize}
For nonlinear models, we use an adaptation of   \cite{Cimen2004} which substitutes the original problem with a finite
sequence of Linear-Quadratic ones to benefit of the points listed
above. 
The linear and nonlinear cases are detailed in the next subsection.

\subsection{Numerical methods for solving (\ref{eq:OC_Cost_Regulator_Form})}
We have defined our state predictor as  solution of a deterministic tracking problem. Several algorithms or tools have been developed to solve this kind of control problem. In this section, $\theta$ is fixed at a given value. Then to simplify the notations, we omit it.

\subsubsection{Linear models \label{sub:linear_SDE_method}}

We consider the case where $\mathbf{A}_{i}:=\mathbf{A}_{\theta,i}=I_{d}+\triangle A_{\theta}(\tilde Z_{i},t_{i}):=I_{d}+\triangle A_{\theta}(t_{i})$
and $\varGamma\left(t_{i}\right):=\varGamma_{\sigma}\left(\tilde Z_{i},t_{i}\right)$ for $i=1, \ldots, n$.
In this framework   $\overline{u}_{Z_{0}}$ is derived  by
solving a finite difference equation known as Riccati equation defined by 
\begin{equation*}
 \begin{array}{lll}
E_{n} & = & C^{T}C\\
h_{n} & = & -C^{T}Y_{n}
\end{array} 
\end{equation*}
and  for $i=n-1, \ldots, 1$: 
\begin{equation}
 \begin{array}{lll}
E_{i}& = & \mathbf{A}_{i}^{T}E_{i+1}\mathbf{A}_{i}+C^{T}C-\triangle \mathbf{A}_{i}^{T}E_{i+1}\varGamma(t_{i})G(E_{i+1})\varGamma(t_{i})^{T}E_{i+1}\mathbf{A}_{i}\\
h_{i} & = & \Delta_{i}\mathbf{A}_{i}^{T}E_{i+1}r(t_{i})+\mathbf{A}_{i}^{T}h_{i+1}-C^{T}Y_{i}\\
 & - & \triangle \mathbf{A}_{i}^{T}E_{i+1}\varGamma(t_{i})G(E_{i+1})\varGamma(t_{i})^{T}\left(h_{i+1}+\Delta_{i}E_{i+1}r(t_{i})\right)\\
\end{array} \label{eq:Riccati_equation}
\end{equation}
with $G(E_{i+1})=\left[\frac{1}{w}\times I_{d}+\triangle \varGamma(t_{i})^{T}E_{i+1}\varGamma(t_{i})\right]^{-1}$. Then  $\overline{u}_{Z_{0},i}$ is given by:
\begin{equation}
\overline{u}_{Z_{0},i}=-\sqrt{\triangle }G(R_{k+1})\varGamma(t_{i})^{T}\left(E_{i+1}\left(\mathbf{A}_{i}\overline{Z}_{i}+\Delta_{i}r(t_{i})\right)+h_{i+1}\right)\label{eq:Optimalcontrol_sequences}
\end{equation}
and $\overline{u}_{Z_{0},i}$ is used in (\ref{eq:sde_original_model})
to obtain $\overline{Z}_{Z_{0}}.$ Moreover,   the initial condition
estimator $\widehat{Z_{0}}$ is   given by $\widehat{Z_{0}}=-\left(E_{0}\right)^{-1}h_{0}.$
Computational details are given in Appendix 2.

\subsubsection{Nonlinear models \label{sub:nonlinear_SDE_method}}
For non-linear model, we propose to apply the previous algorithm thanks to the pseudo-linear representation.  We replace the original problem by a sequence of  Linear-Quadratic control problems solved  iteratively until a convergence criteria is verified. 
This method is an adaptation of  \cite{Cimen2004}
to the case of discrete models. 
We propose the following algorithm to solve   (\ref{eq:OC_Cost_Regulator_Form}). 
\begin{enumerate}
\item Initialisation $\forall i\in\left\llbracket 0,\,n\right\rrbracket ,\,\overline{Z_{i}^{0}}=Z_{0}^{r}$
where $Z_{0}^{r}$ is an arbitrary starting point if $Z_{0}$ is unknown
or $\overline{Z_{0}^{0}}=Z_{0}$ otherwise.
\item At iteration $l$, compute $\left(\overline{Z^{l}},\,\overline{u^{l}}\right)$ by solving the Linear-Quadratic optimal control problem derived from the pseudo-linear representation of the SDE:
\begin{equation}
\begin{array}{ll}
\textrm{Minimize in \ensuremath{u}:} & C_{w}^{l}(u|Y;Z_{0})=\sum_{j=1}^{n-1}\left(\left\Vert C\tilde Z_{i}(Z_{0},u)-Y_{i}\right\Vert _{2}^{2}+\frac{1}{w}u_{i}^{T}u_{i}\right)+\left\Vert C\tilde Z_{n}(Z_{0},u)-Y_{n}\right\Vert _{2}^{2}\\
\textrm{Subject to: } & \left\{ \begin{array}{l}
\tilde Z_{i+1}(Z_{0},u)=\mathbf{A}_{i}^{l}\tilde Z_{i}(Z_{0},u)+\Delta_{i}r(t_{i})+\sqrt{\triangle }\varGamma^{l}(t_{i})u_{i}\\
\tilde Z_{0}=Z_{0}
\end{array}\right.
\end{array}\label{eq:OC_Cost_Regulator_Recursive_Form}
\end{equation}
with $\mathbf{A}_{i}^{l}=I_{d}+\triangle A(\overline{Z_{i}^{l-1}},t_{i})$, $\varGamma^{l}\left(t_{i}\right):=\varGamma\left(\overline{Z_{i}^{l-1}},t_{i}\right)$
and $\overline{Z^{l-1}}$   the state variable corresponding to the
optimal control $\overline{u^{l-1}}$, the minimizer of the cost $C_{w}^{l-1}(u|Y;Z_{0})$.

\item If $\sum_{i=0}^{n}\left\Vert \overline{Z_{i}^{l}}-\overline{Z_{i}^{l-1}}\right\Vert _{2}^{2}<\varepsilon$
then stop otherwise go back to step 2.
\item Use $\overline{Z}\simeq\overline{Z^{l}}$ as state variable predictor.
\end{enumerate}

The interest of this algorithm is that at each iteration $l$, problem (\ref{eq:OC_Cost_Regulator_Recursive_Form})
can be solved  using Linear-Quadratic theory which ensures:
\begin{itemize}
\item the existence and uniqueness of the solution problem $\overline{u^{l}}$, 
\item that $\overline{u^{l}}$ and the corresponding state predictor $\overline{Z^{l}}$
can be computed by solving the Riccati equation (\ref{eq:Riccati_equation}) with
$\mathbf{A}_{i}^{l}:=I_{d}+\triangle A_{\theta}(\overline{Z_{i}^{l-1}},t_{i})$
and $\varGamma^{^{l}}\left(t_{i}\right):=\varGamma\left(\overline{Z_{i}^{l-1}},t_{i}\right)$.
\end{itemize}
  
\subsection{Selection of weights $w$ }
A data driven selection method for $w$ needs to be specified. Let
us denote $\overline{u}_{w}$ the solution of the control problem
(\ref{eq:OC_Cost_Regulator_Form}) obtained for a given weight value
$w$ and the corresponding estimator $\widehat{\psi} = \left(\widehat{\theta},\widehat{\sigma}\right)$.
The sequence $\overline{u}_{w}=(\overline{u}_{w,0},\ldots\overline{u}_{w,n-1})$
is supposed to mimic increments of a Brownian motion. So ideally $\left\Vert \overline{u}_{w,i}\right\Vert _{2}^{2}\sim\chi^{2}(d_{U})$
with $\chi^{2}(d_{U})$ the $\chi^2$ distribution with $d_{U}$ degrees of freedom. 
Thus the i.i.d sequence $(\left\Vert \overline{u}_{w,0}\right\Vert _{2}^{2},\ldots,\left\Vert \overline{u}_{w,n-1}\right\Vert _{2}^{2})$
has ideally a density proportional to $\prod_{i}\left\Vert \overline{u}_{w,i}\right\Vert _{2}^{2\left(\frac{d_{U}}{2}-1\right)}e^{-\left\Vert \overline{u}_{w,i}\right\Vert _{2}^{2}/2}.$
Based on that, we choose the optimal weight $\widehat{w}$ among a
set of values $W$ which maximizes the external criteria:

\begin{equation}
K(w)=\prod_{i}\left\Vert \overline{u}_{w,i}\right\Vert _{2}^{2}{}^{\left(\frac{d_{U}}{2}-1\right)}e^{-\left\Vert \overline{u}_{w,i}\right\Vert _{2}^{2}/2}.
\label{eq:Kw}
\end{equation}
There are some similarities between our method and Generalized Profiling \cite{Ramsay2007}. In both cases, estimation is based on a nested optimization procedure. 1/ The hyperparameter $w$ balancing the model and data fidelity is chosen via the minimization of an outer criteria $K$ given by  (\ref{eq:Kw}). 2/ For a given $w$, the structural parameter $\psi$ is estimated by minimizing the middle criteria $H^{m_{B}}$ given by (\ref{eq:H_mB}). 3/For a given set $(w,\psi)$, a state variable predictor is computed by optimizing the inner criteria $C_w$ given by (\ref{eq:OC_Cost_Regulator_Form}).

\subsection{Summary of  the estimation procedure}

We summarize the method with a pseudo-algorithm formalism: 

\begin{itemize}
\item Outer criteria: estimation of optimal weight $\widehat{w}$ defined
by: 
\[
\widehat{w}:=\arg\min_{w\in W}K(w).
\]

\item Middle criteria: estimate $\widehat{\psi}$ of $\psi$ defined by:
\[
\widehat{\psi}:=\arg\min_{\psi}H^{m_{B}}(\psi\mid Y,\overline{Z}_{\psi}).
\]

\item Inner criteria: compute state predictor $\overline{Z}_{\psi}$ via
algorithm presented in section  \ref{sub:linear_SDE_method} for linear SDEs or section \ref{sub:nonlinear_SDE_method} for nonlinear
ones.
\end{itemize}

\section{Simulation study\label{sec:Simulation}}

\subsection{Experimental design for the simulations}

A simulation study is conducted to analyze the practical accuracy
of our method and its computational cost on three partially observed
hypoelliptic systems, one linear and two non-linear ones. 

Each system is observed on an interval $[0,T]$ and is sampled at
$n$ times uniformly at every $\frac{T}{(n-1)}$ time points with Euler-Maruyama scheme. Three different values for the set  $(T,n)$ are tested to quantify the effects of the interval length $T$ and the sample size $n$ on estimation accuracy. Estimation results are given in terms
of empirical bias and variance computed after Monte- Carlo simulations
based on $N_{MC}=1000$ trials. The sample size has an important impact
on the computational efficiency of the method. Thus, we also give the
mean computational time for a given $w$ .

\subsection{Examples}

\subsubsection{Monotone cyclic feedback system}
We consider a   neuronal   monotone cyclic feedback system proposed in \cite{DitlevsenEva2017}. This model describes the oscillatory behavior of a system of three populations of neurons in interaction:
 \begin{equation}
\begin{array}{l}
dX_1(t)=\left(-\nu X_1(t)+X_2(t)\right)dt\\
dX_2(t)=\left(-\nu X_2(t)+X_3(t)\right)dt\\
dX_3(t)=-\nu X_3(t)dt+cdW_{t}
\end{array}\label{eq:cyclic_neuron_sde}
\end{equation}
where $X_1(t), X_2(t)$ and $X_3(t)$ are the limit dynamics of each population. 

We consider   partial observations with $C=\left(1\,0\,0\right)$.
The model is defined by $g_{\theta}(x_1,x_2,x_3)=\left(\begin{array}{l}
-\nu x_1+x_2\\
-\nu x_2+x_3
\end{array}\right)$ and $h_{\theta}(x_1,x_2,x_3)=-\nu x_3$, 
$A_{\theta}(t)=\left(\begin{array}{ccc}
-\nu & 1 & 0\\
0 & -\nu & 1\\
0 & 0 & -\nu
\end{array}\right),$ $r_{\theta}(t)=\left(0\,0\,0\right)^{T}$ and $\varGamma_{\sigma}(t)=\left(0\,0\,c\right)^{T}.$
The system is   hypoelliptic. Applying the generation step   of the weak H\"{o}rmander condition gives  $L_{0}=\left(0\,0,\,c\right)^{T}$,  $L_{1}=\left[f_{\theta}(Z),L_{0}\right]^{T}=\left(0\,-c\:\nu c\right)^{T}$
and $L_{2}=\left[f_{\theta}(Z),L_{1}\right]^{T}=\left(c\,0\,\nu c^{2}\right)^{T}$. 
Matrix $(L_{0}\,L_{1}\,L_{2})$    spans $\mathbb{R}^{3}.$

Interestingly the need for a second iteration to ensure hypoellipticity
is mimicked by our connexity condition (\ref{eq:connexity_prop})
which needs the auxiliary variable $x_2$ to verify $\frac{\partial g_{\theta,1}}{\partial x_2}\frac{\partial g_{\theta,2}}{\partial x_3}=1\neq0$.

Let us now detail assumption (H1). For $m_{B}=2$,  we have $G_{i+2:i+m_{B}}=0$,
$G_{i+1:i+m_{B}}=0$ and $G_{i:i+m_{B}}=\sqrt{\triangle }c$. Thus
  $m_{B}=2$ implies hypothesis (H1) and  $X_{i:i+m_{B}}$ only depends on $u_{i}$. So the estimator is
defined as the minimizer of $H^{m_{B}}$ given by equation (\ref{eq:H_mB}).

\begin{figure}
\centering{}
\includegraphics[scale=0.7]{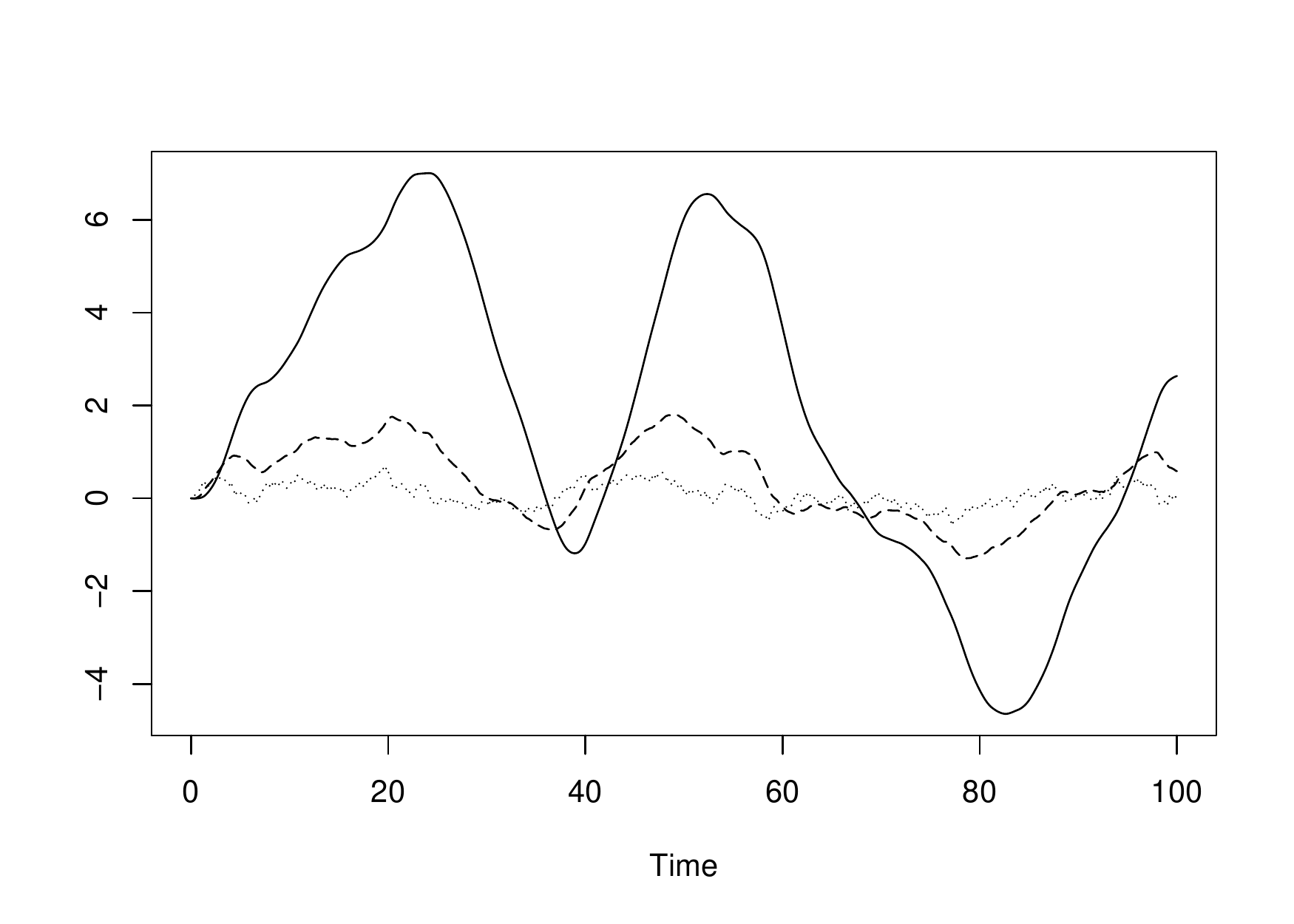}
\caption{\label{fig:cyclic_neuron_sim}Simulated trajectory  of the neuronal   monotone cyclic feedback system (\ref{eq:cyclic_neuron_sde}) with parameters $\nu=0.2$ and $c=0.15$. $X_1:$ Solid line, $X_2:$ dashed line, $X_3:$
dotted line.}
\end{figure}

A thousand simulations are performed with initial conditions  $\left(X_1(0),X_2(0),X_3(0)\right)=\left(0,0,0\right)$ and true parameter values set to $\nu=0.2$ and $c=0.15$. Figure \ref{fig:cyclic_neuron_sim} illustrates a simulation on the observation interval $\left[0,\,T\right]=\left[0,\,20\right].$

As said before, the estimation of $\left(\nu,c\right)$ is made from
the observation of $X_1$ only. The initial condition  are considered known.  We choose $w$ among $W=\left\{ 10^{15},10^{20},\,10^{25},\,10^{30}\right\} .$
Results are given in Table \ref{tab:Cyclic-feedback-system} for different values of $T$ and $n$:  $(T,n)=(10,10^3)$, $(100,10^3)$, $(10,10^4)$. The estimators have 
a good precision which increases with $T$ and mesh refinement. Also,
we point out the computation time is expressed in terms of second
for the case $n=10^{3}$ where it was expressed in terms of hours
in \cite{clairon2020optimal} for an equivalently complex linear model
(Harmonic Oscillator model, section 6.2.1).
\begin{table}
\centering{}
\caption{Cyclic feedback model. Estimation of parameters from 1000 simulated trajectories (mean and variance)  and mean computational time for a given $w$.\label{tab:Cyclic-feedback-system}}
\begin{tabular}{lccc}
\hline 
 & $\nu$ & $c$ & Comp Time \tabularnewline
\hline 
\hline 
True value & 0.2 & 0.15 & \tabularnewline
\hline 
$T=10,n=10^{3}$ & 0.23 (5e-3) & 0.14 (1e-5) & 9s\tabularnewline
\hline 
$T=100,n=10^{3}$ & 0.21 (6e-4) & 0.15 (1e-6) & 19s\tabularnewline
\hline 
$T=10,n=10^{4}$ & 0.22 (2e-3) & 0.15 (1e-6) & 2min50s\tabularnewline
\hline 
\end{tabular}
\end{table}

\subsubsection{Hypoelliptic FitzHugh-Nagumo model}
We consider now the hypoelliptic neuronal model, which is a minimal representation of a spiking neuron model such as the Hodgkin-Huxley model \cite{LeonSamson2017}. It is defined as:
\begin{equation}
\begin{array}{l}
dV_{t}=\frac{1}{\varepsilon}(V_{t}-V_{t}^{3}-U_{t}+s)dt\\
dU_{t}=(\gamma V_{t}-U_{t}+\beta)dt+\sigma dW_{t}
\end{array}\label{eq:hypo_fhn_sde}
\end{equation}
where the variable $V_t$ represents the membrane potential of a neuron at time $t$, and $U_t$ is a recovery variable, which could represent channel kinetics. Parameter $s$ is the magnitude of the stimulus current. Often $s$ represents injected current and is thus controlled in a given experiment. It is therefore assume known and set to $s=0$. 

We consider   partially observations with $C=\left(1\,0\right)$.
 The model is non-linear. Several choices for $A_{\theta}(Z_{t},t)$
are possible, we take $A_{\theta}(Z_{t},t)=\left(\begin{array}{cc}
(1-V_{t}^{2})/\varepsilon & -1/\varepsilon\\
\gamma & -1
\end{array}\right)$ and $r_{\theta}(t)=\left(\begin{array}{c}
0\\
\beta
\end{array}\right)$. Since  $\varGamma_{\sigma}(t)=\left(0\,\sigma\right)^{T}$,
we got $C\varGamma_{\sigma}=0$, $\frac{\partial}{\partial U}(\frac{1}{\varepsilon}(V_{t}-V_{t}^{3}-U_{t}))=-\frac{1}{\epsilon}\neq0$
and $CA_{\theta}(\tilde Z_{i+1},t_{i+1})\varGamma_{\sigma}(\tilde Z_{i},t_{i})=-\sigma/\varepsilon\neq0$.
The model belongs to the case described in Section (\ref{subsec:1_step_hypoelliptic})
and we can use the corresponding simplified expression of $H^{m_{B}}$
for the parameter estimation.  

\begin{figure}
\centering{}
\includegraphics[scale=0.7]{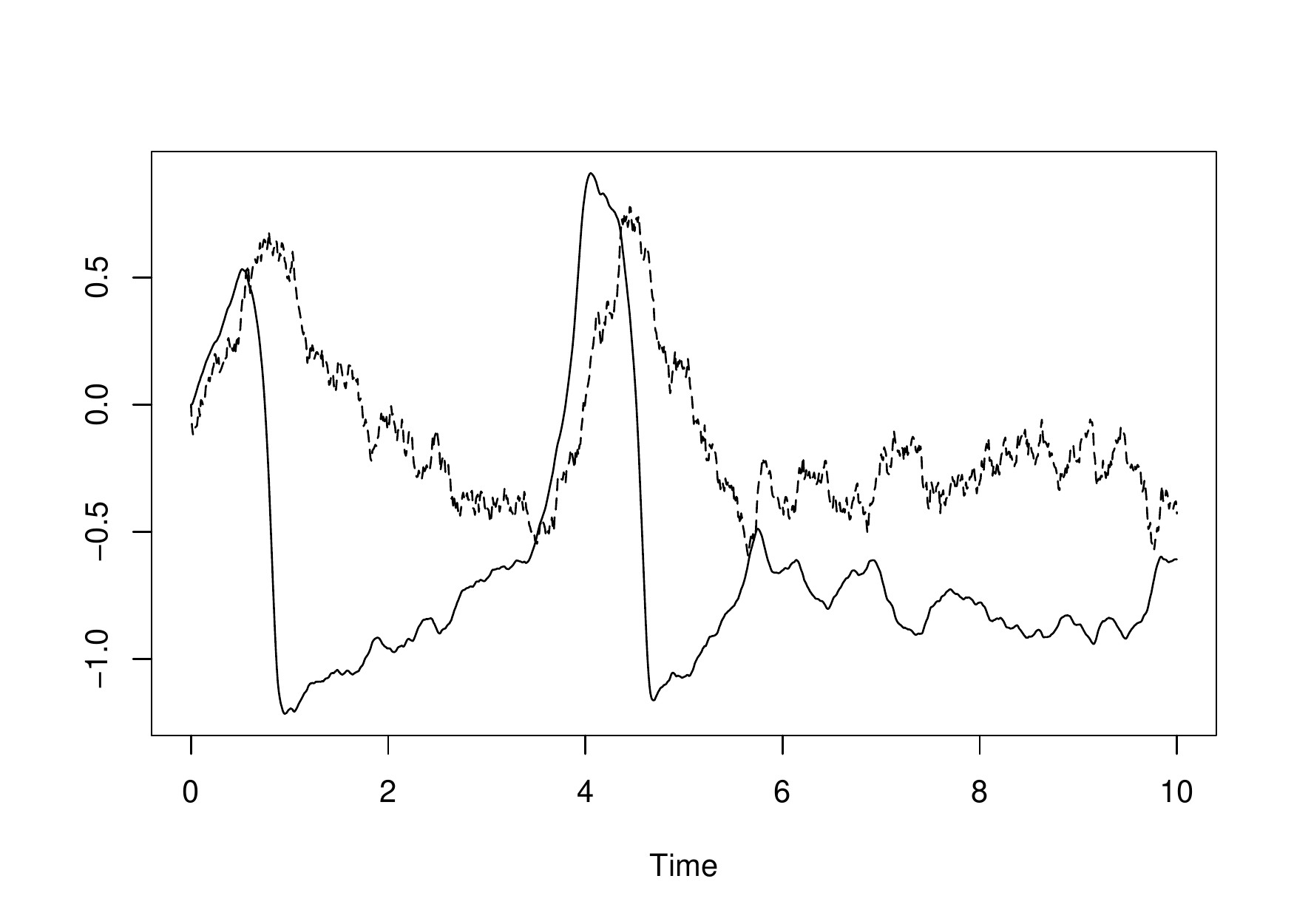}
\caption{\label{fig:fhn_hypo_sim} Simulated trajectory  of the neuronal   hypoelliptic FitzHugh-Nagumo  model (\ref{eq:hypo_fhn_sde}) with parameters $\left(\varepsilon,\gamma,\beta\right)=\left(0.1,1.5,0.8\right)$
and $\sigma=0.3$.
$V:$ solid line, $U:$ dashed line.}
\end{figure}

A thousand simulations are performed with initial conditions set to $\left(V_{0},\,U_{0}\right)=\left(0,0\right)$ and true parameter values $\left(\varepsilon,\gamma,\beta\right)=\left(0.1,1.5,0.8\right)$
and $\sigma=0.3$. Figure \ref{fig:fhn_hypo_sim} illustrates a simulation  on the observation interval $\left[0,\,T\right]=\left[0,\,20\right].$

The estimation is made from the observation of $V$ only and the initial conditions are considered unknown and need to be estimated.  We choose $w$ among $W=\left\{ 10^{16},\,10^{18},\,10^{20},\,10^{25}\right\} .$  

Results are given in Table \ref{tab:Hypoelleptic FitzHugh-Nagumo} for  $(T,n)=(1,10^3), (10,10^3), (1,10^4)$.  As in the previous example, we observe the bias and variance decreasing with $T$ and $n.$ Regarding the computation time and comparative
accuracy with other methods, we recall in Table \ref{tab:Comparison_Hypoelleptic_FitzHugh-Nagumo}, the results obtained for the case $T=10,n=10^{3}$ by  \cite{clairon2020optimal} and  \cite{DitlevsenSamson2019} (they are originally presented in \cite{clairon2020optimal} table 4, section 6.2.2).
We obtain estimation with equivalent or higher accuracy and always
with significantly reduced computational cost.
\begin{table}
\centering{}
\caption{Hypoelliptic FitzHugh-Nagumo model. Estimation of parameters from 1000 simulated trajectories (mean and variance)  and mean computational time for a given $w$.
\label{tab:Hypoelleptic FitzHugh-Nagumo}}
\begin{tabular}{cccccc}
\hline 
 & $\varepsilon$ & $\gamma$ & $\beta$ & $\sigma$ & Comp Time \tabularnewline
\hline 
\hline 
True value & 0.1 & 1.5 & 0.8 & 0.3 & \tabularnewline
\hline 
$T=1,n=10^{3}$ & 0.11 (2e-3) & 2.13 (1e-1) & 1.18 (3e-1) & 0.32 (2e-2) & 1min10s\tabularnewline
\hline 
$T=10,n=10^{3}$ & 0.09 (2e-5) & 1.58 (6e-2) & 0.87 (5e-2) & 0.29 (2e-4) & 3min15s\tabularnewline
\hline 
$T=1,n=10^{4}$ & 0.10 (4e-4) & 1.90 (8e-2) & 1.10 (2e-1) & 0.31 (3e-3) & 8min30s\tabularnewline
\hline 
\end{tabular}
\end{table}

\begin{table}
\centering{}
\caption{Hypoelliptic FitzHugh-Nagumo model. Estimation results for the case $T=10,n=10^{3}$ with our method (first row),   \cite{clairon2020optimal}  (second row), \cite{DitlevsenSamson2019} third row.
\label{tab:Comparison_Hypoelleptic_FitzHugh-Nagumo}}
\begin{tabular}{cccccc}
\hline 
 & $\varepsilon$ & $\gamma$ & $\beta$ & $\sigma$ & Comp Time \tabularnewline
\hline 
\hline 
True value & 0.1 & 1.5 & 0.8 & 0.3 & \tabularnewline
\hline 
$\widehat{\psi}$ & 0.09 (2e-5) & 1.58 (6e-2) & 0.87 (5e-2) & 0.29 (2e-4) & 3min15s\tabularnewline
\hline 
 \cite{clairon2020optimal}   & 0.09 (1e-4) & 1.27 (2e-2) & 0.65 (1e-2) &0.31 (4e-3) & 45min \tabularnewline
\hline 
 \cite{DitlevsenSamson2019}   & 0.10 (1e-4) &  1.59 (2e-2) &0.87 (2e-2) &0.31 (4e-4)  & 6h00min \tabularnewline
\hline 
\end{tabular}
\end{table}

\subsubsection{Synaptic-conductance model}
We consider the conductance-based model with diffusion synaptic input defined in \cite{DitlevsenSamson2019}. It describes the voltage dynamics accross the membrane of a neuron:
\begin{equation}
\begin{array}{l}
C_{c}dV_{t}=(-G_{L}(V_{t}-V_{L})-G_{E,t}(V_{t}-V_{E})-G_{I,t}(V_{t}-V_{I})+I_{inj})dt\\
dG_{E,t}=-\frac{1}{\tau_{E}}(G_{E,t}-g_{E})+\sigma_{E}\sqrt{G_{E,t}}dW_{E,t}\\
dG_{I,t}=-\frac{1}{\tau_{I}}(G_{I,t}-g_{I})+\sigma_{I}\sqrt{G_{I,t}}dW_{I,t}
\end{array}\label{eq:syncon_sde}
\end{equation}
where $C_{c}$ is the total capacitance, $G_L$, $G_E$ and $G_I$ are the leak, excitation and inhibition conductances, $V_L$, $V_E$ and $V_I$ are their respective reversal potentials, and $I_{inj}$ is the injected current. The conductances $G_{E,t}$ and $G_{I,t}$ are assumed to be stochastic functions of time, where $W_{E,t}$ and $W_{I,t}$ are two independent Brownian motions. The square root in the diffusion coefficient ensures that the conductances stay positive. 
Parameters $\tau_E$, $\tau_I$ are time constants, $g_E$, $g_I$ are the mean conductances and $\sigma_E,\sigma_I$ the diffusion coefficients. Here $U_t=(G_{E,t}, G_{I,t})$.  We assume to know the capacitance, the reversal potentials, and the injected current. 

We consider partial observations with only component $V_{t}$ observed which corresponds to the observation matrix $C=\left(1\,0\,0\right)$.  For the pseudo-linear representation, we choose $A_{\theta}(Z_{t},t)=\left(\begin{array}{ccc}
-G_{L}/C_{c} & -(V_{t}-V_{E})/C_{c} & -(V_{t}-V_{I})/C_{c}\\
0 & -1/\tau_{E} & 0\\
0 & 0 & -1/\tau_{I}
\end{array}\right)$,  $r_{\theta}(t)=\left(\begin{array}{c}
\left(G_{L}V_{L}+I_{inj}\right)/C_{c}\\
g_{E}/\tau_{E}\\
g_{I}/\tau_{I}
\end{array}\right).$ Since $g_{\theta}(V,G_{E},G_{I})=(-G_{L}(V-V_{L})-G_{E}(V-V_{E})-G_{I}(V-V_{I})+I_{inj})/C_{c}$,
$C=\left(1\,0\,0\right)$ and $\varGamma_{\sigma}(Z_{t},t)=\left(\begin{array}{ccc}
0 & \sigma_{E}\sqrt{G_{E,t}} & 0\\
0 & 0 & \sigma_{I}\sqrt{G_{I,t}}
\end{array}\right)^{T},$ we got $C\varGamma_{\sigma}=\left(0\,0\right)$, $\frac{\partial g_{\theta}}{\partial\left(G_{E},G_{I}\right)}(V,G_{E},G_{I})$
is of full rank and $CA_{\theta}(\tilde Z_{i+1},t_{i+1})\varGamma_{\sigma}(\tilde Z_{i},t_{i})\neq0$. It enters the framework of section (\ref{subsec:1_step_hypoelliptic}).
 
A thousand simulations are performed with initial conditions set to $\left(V_{0},\,G_{E,0},\,G_{I,0}\right)=\left(-60,\,10,\,1\right)$ and true parameters values set  to $\left(G_{L},V_{L},V_{E},V_{I},I_{inj},g_{E}\right)=\left(50,-70,0,-80,-60,17.8\right)$, $\left(\tau_{E},\,\tau_{I},\,g_{I}\right)=\left(0.5,1,9.4\right)$
and $\left(\sigma_{E},\sigma_{I}\right)=\left(0.1,\,0.1\right)$. The estimation of $\left(\tau_{E},\,\tau_{I},\,g_{I}, \sigma_{E},\sigma_{I}\right)$ is made from the observation of $V$ only. 
For the sake of identifiability, the initial conditions are assumed known.  We choose $w$ among $W=\left\{ 10^{8},5\times10^{8},\,10^{9},\,5\times10^{9}\right\} .$

Results are given in Table \ref{tab:Synaptic-conductance model} for   $(T,n)=(20,10^3), (200,10^3), (20,10^4)$. We observe
a difference in terms of accuracy between the estimation of $\theta=(\tau_{E},\tau_{I},g_{I})$
and $\sigma=\left(\sigma_{E},\sigma_{I}\right)$, the latter suffering
from more bias than the former. This is understandable and already noticed in \cite{DitlevsenSamson2019}. Indeed
we have to estimate the diffusion of two sources of stochastic disturbance
from only one resulting signal. Because of the number of involved
state variables and parameters to estimate, the ratio of observed/unobserved
states and diffusion matrix structure, 
This estimation problem is
more complex than the previous ones as  the number of parameters to estimate is higher and the ratio of observed/unobserved
states is also higher. This explains the higher computational
time for the case $T=20,n=10^{4}$.
\begin{table}
\centering{}
\caption{Synaptic-conductance model. Estimation of parameters from 1000 simulated trajectories (mean and variance)  and mean computational time for a given $w$. \label{tab:Synaptic-conductance model}}
\begin{tabular}{ccccccc}
\hline 
 & $\tau_{E}$ & $\tau_{I}$ & $g_{I}$ & $\sigma_{E}$ & $\sigma_{I}$ & Comp Time\tabularnewline
\hline 
\hline 
True value & 0.5 & 1 & 9.4 & 0.1 & 0.1 & \tabularnewline
\hline 
$T=20,n=10^{3}$ & 0.51 (2e-3) & 1.11 (0.04) & 9.45 (0.02) & 0.05 (3e-4) & 0.15 (1e-3) & 1min40s\tabularnewline
\hline 
$T=200,n=10^{3}$ & 0.47 (4e-4) & 1.11 (0.01) & 9.41 (2e-3) & 0.09 (7e-6) & 0.14 (8e-5) & 2min00s\tabularnewline
\hline 
$T=20,n=10^{4}$ & 0.52 (2e-3) & 1.08 (0.04) & 9.40 (0.02) & 0.11 (4e-5) & 0.07 (1e-3) & 1h40min\tabularnewline
\hline 
\end{tabular}
\end{table}

\section{Discussion}

In this work, we propose a new estimation method which gives an unifying
framework for elliptic and hypoelliptic systems, partially or fully
observed. For this, we rely on a lagged discretization of the original
SDE which let enough time to the stochastic perturbations to affect
all the state variables. By doing so, we have constructed a statistical
criteria based on a well-defined density even for partially observed
hypoelliptic SDEs. This criteria requires a state variable predictor
obtained by solving a control problem balancing data and model fidelity.
The numerical procedures used to solve it do not require $B_{\sigma}(Z_{t},t)$
to be of full rank and so are well adapted to hypoelliptic systems.
Because of this, the deterministic control perspective constitutes
a relevant alternative to MCMC approaches and explain the reasonable
computational cost of our method. It only provides the pointwise estimator
$\overline{Z}$ we need for $\psi$ estimation and does not aim to
reconstruct its whole distribution. Now, we conclude this work by
presenting refinments and extensions of the presented method which
will be investigated in the future.

Our method requires to select an hyperparameter $w$. To bypass this,
we aim to define in future works our predictor $\overline{u}$ as
a maximum a posteriori estimator (MAP) in a functional space. For
this, we will rely on the work of \cite{dashti2013map} in which the
MAP is defined as the solution of a new deterministic optimal control
problem where our regularization term $\frac{2}{w}\ln P(u)$ is replaced
by $\left\Vert u\right\Vert _{E}.$ Here, $E$ denotes the so-called
Cameron-Martin space which is totally determined by the covariance
operator of the Brownian motion $u$ and does not involve a nuisance
parameter. 

Also, we think our method is well suited for generalization to SDEs
driven by processes $W_{t}$ different from the Wiener ones. All it
requires would be to know their densities to modify $H^{m_{B}}$ accordingly
and the penalization term appearing in the control problem. Interestingly
in this general setting, to a given type of SDE will be associated
a given type of deterministic optimal control problem.

We end this work by going one step further about the connection between
our estimation problem and control theory by pointing out its similarity
with the issue of structural assessment of a system controllability
as exposed in \cite{daoutidis1992structural}. Their problem is the
following: given an ODE $\dot{x}=f(x)+Bu$, can we control the behavior
of some system outputs $y_{i}=h_{i}(x)$ knowing that $u=\left\{ u_{j}\right\} _{j\in\left\llbracket 1,d_{u}\right\rrbracket }$
only affect directly a subset of $x$? This leads to define in a similar
way as our $m_{l}$'s the integers $r_{ij}$, named relative orders,
quantifying the sluggishness of $y_{i}$ response to $u_{j}$ variation.
Then, conditions for controllability are formulated via the non-singularity
of a so-called characteristic matrix $C(x)$ constructed from the
$r_{ij}$'s and Lie derivative in a way mirroring the matrix $M$
appearing in the proof of proposition \ref{prop:connexity_prop}.
From this, we hope to derive $m_{B}$ in a less exploratory manner
directly from a given model structure and make explicit its link with
the required number of iteration of generalization step 2 to fulfill H\"{o}rmander condition.

\section*{Code/Software}

Our estimation method is implemented in R and a code reproducing the
examples of Section \ref{sec:Simulation} is available on a GitHub
repository located \href{https://github.com/QuentinClairon/SDE_estimation_via_optimal_control.git}{here}.

\bibliographystyle{Chicago}

\bibliography{biblio_hypo_sde}

\begin{thebibliography}{}

\bibitem[\protect\citeauthoryear{Bierkens, van~der Meulen, and
  Schauer}{Bierkens et~al.}{2020}]{Bierkens2020}
Bierkens, J., F.~van~der Meulen, and M.~Schauer (2020).
\newblock Simulation of elliptic and hypo-elliptic conditional diffusions.
\newblock {\em Adv Applied Proba\/}~{\em 52}.

\bibitem[\protect\citeauthoryear{Buckwar, Samson, Tamborrino, and
  Tubikanec}{Buckwar et~al.}{2021}]{buckwar2021splitting}
Buckwar, E., A.~Samson, M.~Tamborrino, and I.~Tubikanec (2021).
\newblock Splitting methods for sdes with locally lipschitz drift. an
  illustration on the fitzhugh-nagumo model.
\newblock {\em arXiv preprint arXiv:2101.01027\/}.

\bibitem[\protect\citeauthoryear{Cimen}{Cimen}{2008}]{Cimen2008}
Cimen, T. (2008).
\newblock State-dependent riccati equation (sdre) control: A survey.
\newblock {\em IFAC Proceedings\/}~{\em 41}, 3761--3775.

\bibitem[\protect\citeauthoryear{Cimen and Banks}{Cimen and
  Banks}{2004a}]{CimenBanks2004}
Cimen, T. and S.~Banks (2004a).
\newblock Global optimal feedback control for general nonlinear systems with
  nonquadratic performance criteria.
\newblock {\em Systems and Control Letters\/}~{\em 53}, 327--346.

\bibitem[\protect\citeauthoryear{Cimen and Banks}{Cimen and
  Banks}{2004b}]{Cimen2004}
Cimen, T. and S.~Banks (2004b).
\newblock Nonlinear optimal tracking control with application to super-tankers
  for autopilot design.
\newblock {\em Automatica\/}~{\em 40}, 1845--1863.

\bibitem[\protect\citeauthoryear{Clairon}{Clairon}{2021}]{clairon2021regularization}
Clairon, Q. (2021).
\newblock A regularization method for the parameter estimation problem in
  ordinary differential equations via discrete optimal control theory.
\newblock {\em Journal of Statistical Planning and Inference\/}~{\em 210},
  1--19.

\bibitem[\protect\citeauthoryear{Clairon and Brunel}{Clairon and
  Brunel}{2017}]{BrunelClairon_Pontryagin2017}
Clairon, Q. and N.~Brunel (2017).
\newblock Optimal control and additive perturbations help in estimating
  ill-posed and uncertain dynamical systems.
\newblock JASA.

\bibitem[\protect\citeauthoryear{Clairon and Brunel}{Clairon and
  Brunel}{2019}]{clairon2019tracking}
Clairon, Q. and N.~J.-B. Brunel (2019).
\newblock Tracking for parameter and state estimation in possibly misspecified
  partially observed linear ordinary differential equations.
\newblock {\em Journal of Statistical Planning and Inference\/}~{\em 199},
  188--206.

\bibitem[\protect\citeauthoryear{Clairon and Samson}{Clairon and
  Samson}{2020}]{clairon2020optimal}
Clairon, Q. and A.~Samson (2020).
\newblock Optimal control for estimation in partially observed elliptic and
  hypoelliptic linear stochastic differential equations.
\newblock {\em Statistical Inference for Stochastic Processes\/}~{\em 23\/}(1),
  105--127.

\bibitem[\protect\citeauthoryear{Coombes and Byrne}{Coombes and
  Byrne}{2019}]{CoombesByrne2017}
Coombes, S. and A.~Byrne (2019).
\newblock {\em Lecture Notes in Nonlinear Dynamics in Computational
  Neuroscience: from Physics and Biology to ICT}, Chapter Next generation
  neural mass models.
\newblock PoliTO Springer Series. Springer.

\bibitem[\protect\citeauthoryear{Daoutidis and Kravaris}{Daoutidis and
  Kravaris}{1992}]{daoutidis1992structural}
Daoutidis, P. and C.~Kravaris (1992).
\newblock Structural evaluation of control configurations for multivariable
  nonlinear processes.
\newblock {\em Chemical engineering science\/}~{\em 47\/}(5), 1091--1107.

\bibitem[\protect\citeauthoryear{Dashti, Law, Stuart, and Voss}{Dashti
  et~al.}{2013}]{dashti2013map}
Dashti, M., K.~J. Law, A.~M. Stuart, and J.~Voss (2013).
\newblock Map estimators and their consistency in bayesian nonparametric
  inverse problems.
\newblock {\em Inverse Problems\/}~{\em 29\/}(9), 095017.

\bibitem[\protect\citeauthoryear{DeVille, Vanden-Eijnden, and Muratov}{DeVille
  et~al.}{2005}]{DeVille2005}
DeVille, R., E.~Vanden-Eijnden, and C.~Muratov ({2005}, {SEP}).
\newblock {Two distinct mechanisms of coherence in randomly perturbed dynamical
  systems}.
\newblock {\em {Physical Review E}\/}~{\em {72}\/}({3, 1}).

\bibitem[\protect\citeauthoryear{Dietz}{Dietz}{2001}]{Dietz2001}
Dietz, H. (2001).
\newblock Asymptotic behaviour of trajectory fitting estimators for certain
  non-ergodic sde.
\newblock {\em Statistical Inference for Stochastic Processes\/}~{\em 4},
  249--258.

\bibitem[\protect\citeauthoryear{Ditlevsen, Ditlevsen, and Andersen}{Ditlevsen
  et~al.}{2002}]{Ditlevsen2002}
Ditlevsen, P., S.~Ditlevsen, and K.~Andersen (2002).
\newblock The fast climate fluctuations during the stadial and interstadial
  climate states.
\newblock {\em Annals of Glaciology\/}~{\em 35}, 457--462.

\bibitem[\protect\citeauthoryear{Ditlevsen and Greenwood}{Ditlevsen and
  Greenwood}{2013}]{DitlevsenGreenwood2013}
Ditlevsen, S. and P.~Greenwood ({2013}, {AUG}).
\newblock {The Morris-Lecar neuron model embeds a leaky integrate-and-fire
  model}.
\newblock {\em {Journal of Mathematical Biology}\/}~{\em {67}\/}({2}),
  {239--259}.

\bibitem[\protect\citeauthoryear{Ditlevsen and Löcherbach}{Ditlevsen and
  Löcherbach}{2017}]{DitlevsenEva2017}
Ditlevsen, S. and E.~Löcherbach (2017).
\newblock Multi-class oscillating systems of interacting neurons.
\newblock {\em Stochastic Processes and Their Application\/}~{\em 127},
  1840--1869.

\bibitem[\protect\citeauthoryear{Ditlevsen and L{\"o}cherbach}{Ditlevsen and
  L{\"o}cherbach}{2017}]{DitlevsenLocherbach2017}
Ditlevsen, S. and E.~L{\"o}cherbach (2017).
\newblock Multi-class oscillating systems of interacting neurons.
\newblock {\em Stochastic Processes and Their Applications\/}~{\em 127},
  1840--1869.

\bibitem[\protect\citeauthoryear{Ditlevsen and Samson}{Ditlevsen and
  Samson}{2014}]{Ditlevsen2014}
Ditlevsen, S. and A.~Samson (2014).
\newblock Estimation in the partially observed stochastic morris-lecar neuronal
  model with particle filter and stochastic approximation methods.
\newblock {\em Annals of Applied Statistics\/}~{\em 2}, 674--702.

\bibitem[\protect\citeauthoryear{Ditlevsen and Samson}{Ditlevsen and
  Samson}{2019}]{DitlevsenSamson2019}
Ditlevsen, S. and A.~Samson (2019).
\newblock Hypoelliptic diffusions: discretization, filtering and inference from
  complete and partial observations.
\newblock {\em J Royal Statistical Society B\/}~{\em 81}, 361--384.

\bibitem[\protect\citeauthoryear{Ditlevsen and S{\o}rensen}{Ditlevsen and
  S{\o}rensen}{2004}]{Ditlevsen2004}
Ditlevsen, S. and M.~S{\o}rensen (2004).
\newblock Inference for observations of integrated diffusion processes.
\newblock {\em Scand. J. Statist.\/}~{\em 31\/}(3), 417--429.

\bibitem[\protect\citeauthoryear{Engl, Flamm, K{\"u}gler, Lu, M{\"u}ller, and
  Schuster}{Engl et~al.}{2009}]{Engl2009}
Engl, H.~W., C.~Flamm, P.~K{\"u}gler, J.~Lu, S.~M{\"u}ller, and P.~Schuster
  (2009).
\newblock Inverse problems in systems biology.
\newblock {\em Inverse Problems\/}~{\em 25\/}(12).

\bibitem[\protect\citeauthoryear{Genon-Catalot, Jeantheau, and
  Larédo}{Genon-Catalot et~al.}{2000}]{genon-catalot2000}
Genon-Catalot, V., T.~Jeantheau, and C.~Larédo (2000, 12).
\newblock Stochastic volatility models as hidden markov models and statistical
  applications.
\newblock {\em Bernoulli\/}~{\em 6\/}(6), 1051--1079.

\bibitem[\protect\citeauthoryear{Gerstner and Kistler}{Gerstner and
  Kistler}{2002}]{Gerstner2002}
Gerstner, W. and W.~Kistler (2002).
\newblock {\em Spiking Neuron Models}.
\newblock Cambridge University Press.

\bibitem[\protect\citeauthoryear{Gloter}{Gloter}{2000}]{Gloter2000}
Gloter, A. (2000).
\newblock Discrete sampling of an integrated diffusion process and parameter
  estimation of the diffusion coefficient.
\newblock {\em ESAIM Probab. Statist.\/}~{\em 4}, 205--227.

\bibitem[\protect\citeauthoryear{Gloter}{Gloter}{2006}]{Gloter2006}
Gloter, A. (2006).
\newblock Parameter estimation for a discretely observed integrated diffusion
  process.
\newblock {\em Scand. J. Statist.\/}~{\em 33\/}(1), 83--104.

\bibitem[\protect\citeauthoryear{Gloter and Yoshida}{Gloter and
  Yoshida}{2021}]{Gloter2021}
Gloter, A. and N.~Yoshida (2021).
\newblock Adaptive estimation for degenerate diffusion process.
\newblock {\em Electronic J Stat\/}~{\em 15}, 1424--1472.

\bibitem[\protect\citeauthoryear{Goldwyn and Shea-Brown}{Goldwyn and
  Shea-Brown}{2011}]{GoldwynSheaBrown2011}
Goldwyn, J.~H. and E.~Shea-Brown ({2011}, {NOV}).
\newblock {The What and Where of Adding Channel Noise to the Hodgkin-Huxley
  Equations}.
\newblock {\em {PLOS Computational Biology}\/}~{\em {7}\/}({11}).

\bibitem[\protect\citeauthoryear{Graham, Thiery, and Beskos}{Graham
  et~al.}{2019}]{Graham2019}
Graham, M.~M., A.~H. Thiery, and A.~Beskos (2019).
\newblock Manifold markov chain monte carlo methods for bayesian inference in a
  wide class of diffusion models.
\newblock {\em arXiv: Computation\/}.

\bibitem[\protect\citeauthoryear{Iolov, Ditlevsen, and Longtin}{Iolov
  et~al.}{2017}]{Iolov2017}
Iolov, A., S.~Ditlevsen, and A.~Longtin (2017).
\newblock Optimal design for estimation in diffusion processes from first
  hitting times.
\newblock {\em SIAM J. Uncertainty Quantification\/}~{\em 5}, 88--110.

\bibitem[\protect\citeauthoryear{Kutoyants}{Kutoyants}{1991}]{Kutoyants1991}
Kutoyants, Y. (1991).
\newblock Minimum distance parameter estimation for diffusion type observation.
\newblock {\em Comptes rendus de l Academie des sciences\/}~{\em 312}, 637.

\bibitem[\protect\citeauthoryear{Leimkuhler and Matthews}{Leimkuhler and
  Matthews}{2015}]{LeimkuhlerMatthewsBook2015}
Leimkuhler, B. and C.~Matthews (2015).
\newblock {\em Molecular Dynamics with deterministic and stochastic numerical
  methods}, Volume~39 of {\em Interdisciplinary Applied Mathematics}.
\newblock Springer International Publishing Switzerland.

\bibitem[\protect\citeauthoryear{Leon, Rodriguez, and Ruggiero}{Leon
  et~al.}{2019}]{Leon2019}
Leon, J., L.~Rodriguez, and R.~Ruggiero (2019).
\newblock Consistency of a likelihood estimator for stochastic damping
  hamiltonian systems. totally observed data.
\newblock {\em ESAIM PS\/}~{\em 23}, 1--36.

\bibitem[\protect\citeauthoryear{Leon and Samson}{Leon and
  Samson}{2018}]{LeonSamson2017}
Leon, J. and A.~Samson (2018).
\newblock Hypoelliptic stochastic fitzhugh-nagumo neuronal model: mixing,
  up-crossing and estimation of the spike rate.
\newblock {\em Annals of Applied Probability\/}~{\em 28\/}(4), 2243--2274.

\bibitem[\protect\citeauthoryear{Lu, Lin, and Chorin}{Lu et~al.}{2016}]{Lu2016}
Lu, F., K.~Lin, and A.~Chorin (2016).
\newblock comparison of continuous and discrete-time data-based modeling for
  hypoelliptic systems.
\newblock {\em Comm App Math Comp Sci\/}~{\em 11}, 187--216.

\bibitem[\protect\citeauthoryear{Melnykova}{Melnykova}{2020}]{Melnykova2020}
Melnykova, A. (2020).
\newblock Parametric inference for hypoelliptic ergodic diffusions with full
  observations.
\newblock {\em Statistical Inference for Stochastic Processes\/}~{\em 23},
  595--635.

\bibitem[\protect\citeauthoryear{Paninski, Vidne, DePasquale, and
  Fereira}{Paninski et~al.}{2012}]{Paninski2012}
Paninski, L., M.~Vidne, B.~DePasquale, and D.~Fereira ({2012}).
\newblock {Inferring synaptic inputs given a noisy voltage trace via sequential
  Monte Carlo methods}.
\newblock {\em {Journal of Computational Neuroscience}\/}~{\em {33}\/}({1}),
  {1--19}.

\bibitem[\protect\citeauthoryear{Pokern, Stuart, and Wiberg}{Pokern
  et~al.}{2009}]{Pokern2009}
Pokern, Y., A.~Stuart, and P.~Wiberg (2009).
\newblock Parameter estimation for partially observed hypoelliptic diffusions.
\newblock {\em J. Roy. Stat. Soc. B\/}~{\em 71\/}(1), 49--73.

\bibitem[\protect\citeauthoryear{Ramsay, Hooker, Cao, and Campbell}{Ramsay
  et~al.}{2007}]{Ramsay2007}
Ramsay, J., G.~Hooker, J.~Cao, and D.~Campbell (2007).
\newblock Parameter estimation for differential equations: A generalized
  smoothing approach.
\newblock {\em Journal of the Royal Statistical Society (B)\/}~{\em 69},
  741--796.

\bibitem[\protect\citeauthoryear{Samson and Thieullen}{Samson and
  Thieullen}{2012}]{Samson2012}
Samson, A. and M.~Thieullen (2012).
\newblock Contrast estimator for completely or partially observed hypoelliptic
  diffusion.
\newblock {\em Stochastic Processes and Their Applications\/}~{\em 122},
  2521--2552.

\bibitem[\protect\citeauthoryear{Tuckwell and Ditlevsen}{Tuckwell and
  Ditlevsen}{2016}]{TuckwellDitlevsen2016}
Tuckwell, H.~C. and S.~Ditlevsen ({2016}, {OCT}).
\newblock {The Space-Clamped Hodgkin-Huxley System with Random Synaptic Input:
  Inhibition of Spiking by Weak Noise and Analysis with Moment Equations}.
\newblock {\em {Neural Computation}\/}~{\em {28}\/}({10}), {2129--2161}.

\bibitem[\protect\citeauthoryear{Wu}{Wu}{2001}]{Wu2001}
Wu, L. (2001).
\newblock Large and moderate deviations and exponential convergence for
  stochastic damping hamiltonian systems.
\newblock {\em Stochastic Process. Appl.\/}~{\em 91}, 205--238.

\end{thebibliography}


\begin{thebibliography}{}

\bibitem[\protect\citeauthoryear{Bertsekas}{Bertsekas}{2005}]{Bertsekas2005}
Bertsekas, D. (2005).
\newblock {\em Dynamic Programming and Optimal Control}.
\newblock Athena Scientific.

\end{thebibliography}
\end{document}


\def\spacingset#1{\renewcommand{\baselinestretch}%
{#1}\small\normalsize} \spacingset{1}


\if1\blind
{
  \title{\bf  Supplementary Materials for "Optimal control   for parameter estimation in partially observed  hypoelliptic stochastic differential equations"}
  \author{Quentin Clairon  
\hspace{.2cm}\\
  University of Bordeaux, Inria Bordeaux Sud-Ouest, \\
Inserm, Bordeaux Population Health Research Center, SISTM Team, UMR1219, \\
F-33000 Bordeaux, France\\
Vaccine Research Institute, F-94000 Créteil, France\\
    and \\
    Adeline Samson  \thanks{
    The author gratefully acknowledges the LabEx PERSYVAL-Lab
(ANR-11-LABX-0025-01) and  MIAI@Grenoble Alpes, (ANR-19-P3IA-0003).}\\
    Univ. Grenoble Alpes, CNRS, Inria, Grenoble  INP*,  LJK,\\
  38000 Grenoble, France\\
* Institute of Engineering Univ. Grenoble Alpes
}
  \maketitle
} \fi

\if0\blind
{
  \bigskip
  \bigskip
  \bigskip
  \begin{center}
    {\LARGE\bf  Supplementary Materials for "Optimal control   for parameter estimation in partially observed  hypoelliptic stochastic differential equations"}
\end{center}
  \medskip
} \fi

\bigskip

\newpage
\spacingset{1.45} 

\section*{Appendix 1: proof of propositions 1 and 2}

In this section we give the proof of propositions 1 and 2 on which we build our lagged
constrast function $H^{m_{B}}$. For the sake of clarity we restrict
ourselves in a region $\mathcal{Z}\subseteq\mathbb{R}^{d}$ 

for each
variable $j\in\left\llbracket 1,\,d\right\rrbracket $ we got  either:
\begin{enumerate}
\item  $\forall Z\in\mathcal{Z}$, $\frac{\partial g_{\theta}(Z)}{\partial Z_{j}}=0$ 

or
\item $\forall Z\in\mathcal{Z}$, $\frac{\partial g_{\theta}(Z)}{\partial Z_{j}}\neq0$.
\end{enumerate}
\begin{prop*}
If the SDE (1) is hypoelliptic then for
each $l\in\left\llbracket 1,d_{V}\right\rrbracket $ it exists $j\in\left\llbracket 1,d_{U}\right\rrbracket $
and a sequence $\left(q_{1},\ldots q_{m_{l}}\right)\in\left\llbracket 1,\,d_{V}\right\rrbracket ^{m_{l}}$
with $m_{l}\leq d_{V}$ such that 
\[
\frac{\partial g_{\theta,l}(Z_{t},t)}{\partial V_{q_{1}}}\frac{\partial g_{\theta,q_{1}}(Z_{t},t)}{\partial V_{q_{2}}}\cdots\frac{\partial g_{q_{\theta,m_{l}-1}}(Z_{t},t)}{\partial V_{q_{m_{l}}}}\frac{\partial g_{\theta,q_{m_{l}}}(Z_{t},t)}{\partial U_{j}}\neq0.
\]
\end{prop*}
\begin{proof}
In this proof, we drop the dependence in $t$ for the sake of clarity.
Since $B_{\sigma}(Z_{t})$ is non-singular for every possible values
$Z_{t}\in\mathcal{Z}$, then the vectors $\varGamma_{j}$ already
span $\mathbb{R}^{d_{U}}.$ Therefore hypoellipticity is ensured as
soon as there are vectors belonging to $\mathcal{L}$ which have their
first $d_{V}$ entries spanning $\mathbb{R}^{d_{V}}.$ 

Let us construct recursively the vectors defined by the generation
rule of the H\"{o}rmander condition, we denote $P^{(k)}$ a 
vector constructed at the $k$-th iteration. They are defined by either
$P^{(k)}=\left[f(Z_{t}),\,P^{(k-1)}\right]$ or $P^{(k)}=\left[\varGamma_{j}(Z_{t}),\,P^{(k-1)}\right]$
for $1\leq j\leq d_{U}$ with $\left[.,.\right]$ the Lie bracket. Each $P^{(k)}$ is a $d$ vector composed
of two subvectors $P^{(k)}=\left(\begin{array}{l}
P_{V}^{(k)}\\
P_{U}^{(k)}
\end{array}\right)$ and we denote $\left\{ P_{V}^{(k)}\right\} $ the set of all generated
subvectors at iteration $k$. Hypoellipticity implies that there is
an index $k$ such that the recursively defined matrix $M^{(k)}=\left[\left\{ P_{V}^{(k)}\right\} ,\,M^{\left(k-1\right)}\right]$
is of rank $d_{V}$ and especially the matrix $M^{(k)}$ has only
non null rows.

We want to prove the proposition for a fixed index $l$. Let $P_{l,V}^{(k)}$
denote the $l-th$ entry of $P_{V}^{(k)}$ which is located in the
$l-$th row of $M^{(k)}$. We now prove recursively the proposition
$\mathcal{P}(k)$: if the sequence $\left\{ P_{l,V}^{(k)}\right\} _{k\in\mathbb{N}}$
becomes non null at the $k+1$-th iteration i.e if there is $(k+1)$
such that $P_{l,V}^{(k+1)}\neq0$ and $P_{l,V}^{(k)}=0$ for all $k$,
then it exists $j\in\left\llbracket 1,d_{U}\right\rrbracket $ and
a sequence $\left(q_{1},\ldots q_{m_{l}}\right)$ with $m_{l}\leq k$
such that $\frac{\partial g_{l}(Z)}{\partial V_{q_{1}}}\frac{\partial g_{q_{1}}(Z)}{\partial V_{q_{2}}}\cdots\frac{\partial g_{q_{m_{l}}}(Z)}{\partial U_{j}}\neq0.$

Initialization: for $k=0$, $m_{l}=0$ and we only have to prove the
existence of $j$ such that $\frac{\partial g_{l}(Z)}{\partial U_{j}}\neq0.$
The only vectors in $\mathcal{L}$ before any application of the generation
step are the columns $\varGamma_{h}$ with $h=1,\ldots d_{U}$. From
this, there are 1) $d_{U}^{2}$ vectors $P^{(1)}$ given by $P^{(1)}=\left[\varGamma_{h},\,\varGamma_{j}\right]$
and 2) $d_{U}$ vectors given by $P^{(1)}=\left(\begin{array}{l}
P_{V}^{(1)}\\
P_{U}^{(1)}
\end{array}\right)=\left[f(Z),\,\varGamma_{h}\right].$ In case 1), the $d_{V}$ first entries are given by $\left[\varGamma_{h},\,\varGamma_{j}\right]_{r}=\left(\frac{\partial\varGamma_{r,j}(Z)}{\partial z}\right)\varGamma_{h}-\left(\frac{\partial\varGamma_{r,h}(Z)}{\partial z}\right)\varGamma_{j}=0$
for $r\in\left\llbracket 1,\,d_{V}\right\rrbracket $ and so $P_{V}^{(1)}=0$.
So let us focus on case 2) to find non-zero $P_{V}^{(1)}$. In case
2), we got $P_{l,V}^{(1)}=\left[f(Z),\,\varGamma_{h}\right]_{l}=-\frac{\partial g_{l}(Z)}{\partial z}^{T}\left(\begin{array}{c}
0_{d_{V},1}\\
B_{h}(Z)
\end{array}\right)=-\frac{\partial g_{l}(Z)}{\partial U}^{T}B_{h}(Z)=-\sum_{j}\frac{\partial g_{l}(Z)}{\partial U_{j}}B_{j,h}(Z).$ So if $P_{l,V}^{(1)}\neq0$, it imposes the existence of at least
one $j$ and $h$ such that $\frac{\partial g_{l}(Z)}{\partial U_{j}}B_{j,h}(Z)\neq0$
and so $\frac{\partial g_{l}(Z)}{\partial U_{j}}\neq0$ which concludes
the initialization.

Recursion: let us assume that $\mathcal{P}$ is true up to iteration
$k-1$. Here, we consider that $\frac{\partial g_{l}(Z)}{\partial U_{j}}=0$
for all $j\in\left\llbracket 1,d_{U}\right\rrbracket $ otherwise
the proof is trivial. Again, we consider two cases either 1)$P^{(k)}=\left[\varGamma_{j}(Z_{t}),\,P^{(k-1)}\right]$
or 2)$P^{(k)}=\left[f(Z_{t}),\,P^{(k-1)}\right].$ We start with the
case 1) with vectors $P^{(k)}=\left[\varGamma_{j}(Z_{t}),\,P^{(k-1)}\right]$,
now let us denote $r\in\left\llbracket 1,\,d_{V}\right\rrbracket $
such that $P_{r}^{(k-1)}=0$. Then $P_{r}^{(k)}$, the $r-$th entry
of $P_{r}^{(k)}$ is given by:
\[
P_{r}^{(k)}=\left[\varGamma_{j}(Z),\,P^{(k-1)}\right]_{r}=\left(\frac{\partial P_{r}^{(k-1)}}{\partial z}\right)^{T}\left(\begin{array}{c}
0_{d_{V},1}\\
B_{j}(Z)
\end{array}\right)-\left(\frac{\partial\varGamma_{r,j}(Z)}{\partial z}\right)^{T}P^{(k-1)}=0.
\]
 Thus $P^{(k)}$ has the same null entries as $P^{(k-1)}.$ Now let
us focus on case 2), let us consider $P_{l,V}^{(k+1)}$ such that
$P_{l,V}^{(k)}=0$. The recursive definition of $P_{l,V}^{(k+1)}$
gives:
\[
\begin{array}{lll}
P_{l,V}^{(k+1)} & = & \left(\frac{\partial P_{l,V}^{(k)}}{\partial z}\right)^{T}f(Z)-\frac{\partial g_{l}(Z)^{T}}{\partial V}P_{V}^{(k)}-\frac{\partial g_{l}(Z)^{T}}{\partial U}P_{U}^{(k)}\\
 & = & -\sum_{p}\frac{\partial g_{l}(Z)}{\partial V_{p}}P_{p,V}^{(k)}.
\end{array}
\]
 Thus $P_{l,V}^{(k+1)}\neq0$ imposes the existence of $p$ such that
$\frac{\partial g_{l}(Z)}{\partial V_{p}}\neq0$ and $P_{p,V}^{(k)}\neq0$.
Let us denote $k_{1}<k$ the first index such that $P_{p,V}^{(k_{1})}=0$
as well as $P_{p,V}^{(k_{2}+1)}\neq0$ for $k_{1}\leq k_{2}<k.$ By
using the recursion hypothesis, we derive the existence of a sequence
$\left(q_{1},\ldots q_{m_{p}}\right)$ with $m_{p}\leq k_{1}$ such
that $\frac{\partial g_{p}(Z)}{\partial V_{q_{1}}}\frac{\partial g_{q_{1}}(Z)}{\partial V_{q_{2}}}\cdots\frac{\partial g_{q_{m_{p}-1}}(Z)}{\partial V_{q_{m_{p}}}}\frac{\partial g_{q_{m_{p}}}(Z)}{\partial U_{j}}\neq0$.
By multiplying the last quantity with $\frac{\partial g_{l}(Z)}{\partial V_{p}}$,
the sequence $\left(p,q_{1},\ldots q_{m_{p}}\right)$ of length less
thant $k$ is such that $\frac{\partial g_{l}(Z)}{\partial V_{p}}\frac{\partial g_{p}(Z)}{\partial V_{q_{1}}}\frac{\partial g_{q_{1}}(Z)}{\partial V_{q_{2}}}\cdots\frac{\partial g_{q_{m_{p}-1}}(Z)}{\partial V_{q_{m_{p}}}}\frac{\partial g_{q_{m_{p}}}(Z)}{\partial U_{j}}\neq0$
which completes the proof for the recursion.
\end{proof}
We demonstrated the existence of paths between each smooth variable
to the rough ones in the continuous setting. Now, we derive the consequence
of this existence on the discretized model used for   parameter
estimation.

\begin{prop*}
 Let us assume that for each
$l\in\left\llbracket 1,d_{_{V}}\right\rrbracket $ it exists a unique
minimal length sequence $\left(q_{1},\ldots q_{m_{l}}\right)$ such
that the connexity property holds. Then
the Euler-Maruyama approximation ($\tilde Z_{i}$, for $i=1, \ldots, n$)  is such that for each $l\in\left\llbracket 1,d_{_{V}}\right\rrbracket$ it exists $j\in\left\llbracket 1,d_{U}\right\rrbracket $ such that:

 $$\frac{\partial \tilde V_{l,i+k}}{\partial \tilde U_{j,i}}=0, \quad \mbox{for } k=1, \ldots, m_l \quad \mbox{and}\quad \frac{\partial \tilde V_{l,i+m_{l}+1}}{\partial \tilde U_{j,i}}\neq0$$
\end{prop*}
\begin{proof}
Again for the sake of notational simplicity, we drop the dependence
in $\theta$ and $t$ in drift functions $g_{\theta}$ and $h_{\theta}$. The proof is based on the general relationship:
\[
\begin{array}{lll}
\frac{\partial \tilde V_{l,i+r+1}}{\partial \tilde U_{j,i}} & = & \frac{\partial \tilde V_{l,i+r}}{\partial \tilde U_{j,i}}+\triangle\left(\sum_{p}\frac{\partial g_{l}}{\partial V_{p}}(\tilde Z_{i+r})\frac{\partial \tilde V_{p,i+r}}{\partial \tilde U_{j,i}}+\sum_{p_{U}}\frac{\partial g_{l}}{\partial U_{p_{U}}}(\tilde Z_{i+r})\frac{\partial \tilde U_{p_{U},i+r}}{\partial \tilde U_{j,i}}\right)\\
\frac{\partial \tilde U_{l_{U},i+r+1}}{\partial \tilde U_{j,i}} & = & \frac{\partial \tilde U_{l_{U},i+r}}{\partial \tilde U_{j,i}}+\triangle\left(\sum_{p}\frac{\partial h_{l_{U}}}{\partial  V_{p}}(\tilde Z_{i+r})\frac{\partial \tilde V_{p,i+r}}{\partial \tilde U_{j,i}}+\sum_{p_{U}}\frac{\partial h_{l_{U}}}{\partial U_{p_{U}}}(\tilde Z_{i+r})\frac{\partial \tilde U_{p_{U},i+r}}{\partial \tilde U_{j,i}}\right)\\
 & + & \sqrt{\triangle}\frac{\partial B_{l_{U}}}{\partial U_{j,i}}(\tilde Z_{i+r})u_{i+r}
\end{array}
\]
obtained by differentiation of the discretized SDE for $l\in\left\llbracket 1,\,d_{V}\right\rrbracket $
and $l_{U}\in\left\llbracket 1,\,d_{U}\right\rrbracket $. Since the
discretized variables can not be instantaneously affected by other
ones, we get the initial values: $\frac{\partial \tilde V_{p,i}}{\partial \tilde U_{j,i}}=0$,
$\frac{\partial \tilde U_{p_{U},i}}{\partial \tilde U_{j,i}}=0$ and $\frac{\partial \tilde V_{j,i}}{\partial \tilde V_{j,i}}=1$
for $p\in\left\llbracket 1,\,d_{V}\right\rrbracket $ and $p_{U}\in\left\llbracket 1,\,d_{U}\right\rrbracket \setminus\left\{ j\right\} $. Then
we derive:
\[
\begin{array}{lll}
\frac{\partial \tilde V_{l,i+1}}{\partial \tilde U_{j,i}} & = & \triangle\left(\sum_{p}\frac{\partial g_{l}}{\partial V_{p}}(\tilde Z_{i})\frac{\partial \tilde V_{p,i}}{\partial \tilde U_{j,i}}+\sum_{p_{U}}\frac{\partial g_{l}}{\partial U_{p_{U}}}(\tilde Z_{i})\frac{\partial \tilde U_{p_{U},i}}{\partial \tilde U_{j,i}}\right)=\triangle\frac{\partial g_{l}}{\partial U_{j}}(\tilde Z_{i})\\
\frac{\partial \tilde U_{l_{U},i+1}}{\partial \tilde U_{j,i}} & = & \boldsymbol{1}_{\left\{ l_{U}=j\right\} }+\triangle\left(\sum_{p}\frac{\partial h_{l_{U}}}{\partial  V_{p}}(\tilde Z_{i})\frac{\partial \tilde V_{p,i}}{\partial \tilde U_{j,i}}+\sum_{p_{U}}\frac{\partial h_{l_{U}}}{\partial U_{p_{U}}}(\tilde Z_{i})\frac{\partial \tilde U_{p_{U},i}}{\partial \tilde U_{j,i}}\right)+\sqrt{\triangle}\frac{\partial B_{l_{V}}}{\partial U_{j}}(\tilde Z_{i})u_{i}\\
 & = & \boldsymbol{1}_{\left\{ l_{U}=j\right\} }+\triangle\frac{\partial h_{l_{U}}}{\partial U_{j}}(\tilde Z_{i})+\sqrt{\triangle}\frac{\partial B_{l_{V}}}{\partial  U_{j}}(\tilde Z_{i})u_{i}.
\end{array}
\]
Based on the previous equalities, we can prove our claim for the special
case $m_{l}=0$, we have $\frac{\partial \tilde V_{l,i+m_{l}+1}}{\partial \tilde U_{j,i}}=\frac{\partial \tilde V_{l,i+1}}{\partial \tilde U_{j,i}}\neq0$
if and only if $\frac{\partial g_{l}}{\partial U_{j}}\neq0$. For
the general case $m_{l}>0$, we prove by recursion the following proposition
$\mathcal{P}(m)$: `` for each $l\in\left\llbracket 1,d_{V}\right\rrbracket $
if it exists an unique minimal length sequence $\left(q_{1}^{l},\ldots q_{m_{l}}^{l}\right)$
of length $m_{l}\leq m$ such that $\frac{\partial g_{l}(Z)}{\partial V_{q_{1}^{l}}}\cdots\frac{\partial g_{q_{m_{l}}^{l}}(Z)}{\partial U_{j}}\neq0$,
then $\frac{\partial \tilde V_{l,i+k}}{\partial  \tilde U_{j,i}}=0$ for $k\leq m_l$
and $\frac{\partial  \tilde V_{l,i+m_{l}+1}}{\partial  \tilde U_{j,i}}\neq0$''

Let us prove $\mathcal{P}(m)$ for $m=1$ for the initialization.
Based on the general relationship for $\frac{\partial  \tilde V_{l,i+r+1}}{\partial  \tilde U_{j,i}}$
and the results we got for $\frac{\partial  \tilde V_{l,i+1}}{\partial  \tilde U_{j,i}}$
and $\frac{\partial  \tilde U_{l_{V},i+1}}{\partial  \tilde U_{j,i}}$, we derive:
\[
\begin{array}{lll}
\frac{\partial  \tilde V_{l,i+2}}{\partial  \tilde U_{j,i}} & = & \frac{\partial  \tilde V_{l,i+1}}{\partial  \tilde U_{j,i}}+\triangle\left(\sum_{p}\frac{\partial g_{l}}{\partial  V_{p}}(  \tilde Z_{i+1})\frac{\partial  \tilde V_{p,i+1}}{\partial  \tilde U_{j,i}}+\sum_{p_{U}}\frac{\partial g_{l}}{\partial  U_{p_{U}}}( \tilde Z_{i+1})\frac{\partial  \tilde U_{p_{U},i+1}}{\partial  \tilde U_{j,i}}\right)\\
 & = & \triangle\frac{\partial g_{l}}{\partial  U_{j}}( \tilde Z_{i})+\triangle^{2}\sum_{p}\frac{\partial g_{l}}{\partial  V_{p}}( \tilde Z_{i+1})\frac{\partial g_{p}}{\partial  U_{j}}( \tilde Z_{i})\\
 & + & \triangle\sum_{p_{U}}\frac{\partial g_{l}}{\partial  U_{p_{U}}}( \tilde Z_{i+1})\frac{\partial  \tilde U_{p_{U},i+1}}{\partial  \tilde U_{j,i}}
\end{array}
\]
Since $m_{l}>0$ and $m=1$, the length of the minimal sequence
$m_{l}$ is necessary $1$. So there is no direct link between $ \tilde V_{l}$
and the rough part. That means $\frac{\partial g_{l}}{\partial U_{j}}( \tilde Z_{i})=0$
and $\frac{\partial g_{l}}{\partial  U_{p_{U}}}( \tilde Z_{i+1})=0$ and
we obtain $\frac{\partial  \tilde V_{l,i+1}}{\partial  \tilde U_{j,i}}=\triangle\frac{\partial g_{l}}{\partial  V_{j}}( \tilde Z_{i})=0$
as well as the simplified expression $\frac{\partial  \tilde V_{l,i+2}}{\partial  \tilde U_{j,i}}=\triangle^{2}\sum_{p}\frac{\partial g_{l}}{\partial V_{p}}( \tilde Z_{i+1})\frac{\partial g_{p}}{\partial  U_{j}}( \tilde Z_{i}).$
By assumption, the minimal sequence is unique, this imposes $\frac{\partial  \tilde U_{l,i+2}}{\partial  \tilde V_{j,i}}=\triangle^{2}\sum_{p}\frac{\partial g_{l}}{\partial  V_{p}}( \tilde Z_{i+1})\frac{\partial g_{p}}{\partial U_{j}}( \tilde Z_{i})=\triangle^{2}\frac{\partial g_{l}}{\partial  V_{q_{1}^{l}}}( \tilde Z_{i+1})\frac{\partial g_{q_{1}^{l}}}{\partial  U_{j}}( \tilde Z_{i})\neq0$
which ends the initialization. Now let us assume $\mathcal{P}(m)$
is true. Let us choose an arbitrary $l\in\left\llbracket 1,d_{_{U}}\right\rrbracket $
and let us assume there is a unique sequence $\left(q_{1}^{l},\ldots q_{m_{l}+1}^{l}\right)$
of minimal length $m_{l}+1$ such that $\frac{\partial g_{l}(Z)}{\partial V_{q_{1}^{l}}}\cdots\frac{\partial g_{q_{m_{l+1}}^{l}}(Z)}{\partial U_{j}}\neq0$.
The general relationship gives us:

\[
\begin{array}{l}
\frac{\partial  \tilde V_{l,i+m_{l}+2}}{\partial  \tilde U_{j,i}}=\frac{\partial  \tilde V_{l,i+m_{l}+1}}{\partial  \tilde U_{j,i}}+\triangle\left(\sum_{p}\frac{\partial g_{l}}{\partial  V_{p}}(Z_{i+m_{l}+1})\frac{\partial  \tilde V_{p,i+m_{l}+1}}{\partial  \tilde U_{j,i}}+\sum_{p_{U}}\frac{\partial g_{l}}{\partial  U_{p_{U}}}( \tilde Z_{i+m_{l}+1})\frac{\partial  \tilde U_{p_{U},i+m_{l}+1}}{\partial  \tilde U_{j,i}}\right)\\
\frac{\partial  \tilde V_{l,i+m_{l}+1}}{\partial  \tilde U_{j,i}}=\frac{\partial  \tilde V_{l,i+m_{l}}}{\partial  \tilde U_{j,i}}+\triangle\left(\sum_{p}\frac{\partial g_{l}}{\partial V_{p}}( \tilde Z_{i+m_{l}})\frac{\partial \tilde V_{p,i+m_{l}}}{\partial  \tilde U_{j,i}}+\sum_{p_{U}}\frac{\partial g_{l}}{\partial  U_{p_{U}}}( \tilde Z_{i+m_{l}})\frac{\partial  \tilde U_{p_{U},i+m_{l}}}{\partial  \tilde U_{j,i}}\right)
\end{array}
\]
Necessarily $\frac{\partial g_{l}}{\partial  U_{p_{U}}}=0$ and by
recursion hypothesis, we got  $\frac{\partial  \tilde V_{l,i+k}}{\partial  \tilde U_{j,i}}=0$ for every $k\leq m_l$ since there is no minimal sequence of length $m_l$ or lower,
so:
\[
\begin{array}{l}
\frac{\partial  \tilde V_{l,i+m_{l}+2}}{\partial  \tilde U_{j,i}}=\frac{\partial  \tilde V_{l,i+m_{l}+1}}{\partial  \tilde U_{j,i}}+\triangle\sum_{p}\frac{\partial g_{l}}{\partial  V_{p}}( \tilde Z_{i+m_{l}+1})\frac{\partial  \tilde V_{p,i+m_{l}+1}}{\partial  \tilde U_{j,i}}\\
\frac{\partial  \tilde V_{l,i+m_{l}+1}}{\partial  \tilde U_{j,i}}=\triangle\sum_{p}\frac{\partial g_{l}}{\partial V_{p}}( \tilde Z_{i+m_{l}})\frac{\partial  \tilde V_{p,i+m_{l}}}{\partial  \tilde U_{j,i}}.
\end{array}
\]
Now let us assume $\frac{\partial  \tilde V_{l,i+m_{l}+1}}{\partial  \tilde U_{j,i}}\neq0$,
this implies there is at least one $p$ such that $\frac{\partial g_{l}}{\partial V_{p}}(Z)\neq0$
and $\frac{\partial  \tilde V_{p,i+m_{l}}}{\partial  \tilde U_{j,i}}\neq0$. By
contraposition of the recursion hypothesis if there is no minimal
sequence $\left(q_{1}^{p},\ldots q_{m_{l}-1}^{p}\right)$ such that
$\frac{\partial g_{p}(Z)}{\partial V_{q_{1}^{p}}}\cdots\frac{\partial g_{q_{m_{l}-1}^{p}}(Z)}{\partial U_{j}}\neq0$,
then $\frac{\partial  \tilde V_{p,i+m_{l}}}{\partial  \tilde U_{j,i}}=0$. So necessarily
there is a sequence such that $\frac{\partial g_{p}(Z)}{\partial V_{q_{1}^{p}}}\cdots\frac{\partial g_{q_{m_{l}-1}^{p}}(Z)}{\partial U_{j}}\neq0$
to ensure $\frac{\partial V_{p,i+m_{l}}}{\partial U_{j,i}}\neq0$.
By multiplying this last sequence by $\frac{\partial g_{l}}{\partial V_{p}}(Z)$,
we conclude   the existence of the sequence $\left(p,q_{1}^{p},\ldots q_{m_{l}-1}^{p}\right)$
of length $m_{l}$ such that $\frac{\partial g_{l}}{\partial V_{p}}(Z)\frac{\partial g_{p}(Z)}{\partial V_{q_{1}^{p}}}\cdots\frac{\partial g_{q_{m_{l}-1}^{p}}(Z)}{\partial U_{j}}\neq0$
in contradiction with the minimality of $\left(q_{1}^{l},\ldots q_{m_{l}+1}^{l}\right)$. So by contradiction $\frac{\partial  \tilde V_{l,i+m_{l}+1}}{\partial  \tilde U_{j,i}}=0$.
From this, we get: $\frac{\partial  \tilde V_{l,i+m_{l}+2}}{\partial  \tilde U_{j,i}}=\triangle\sum_{p}\frac{\partial g_{l}}{\partial  V_{p}}( \tilde Z_{i+m_{l}})\frac{\partial  \tilde V_{p,i+m_{l}}}{\partial  \tilde U_{j,i}}$,
assuming $\frac{\partial  \tilde V_{l,i+m_{l}+2}}{\partial  \tilde U_{j,i}}=0$ leads
to have for all $p\in\left\llbracket 1,\,d_{V}\right\rrbracket $
and sequences $\left(q_{1}^{p},\ldots q_{m_{l}}^{p}\right)$ the equality
$\frac{\partial g_{p}(Z)}{\partial V_{q_{1}^{p}}}\cdots\frac{\partial g_{q_{m_{l}}^{p}}(Z)}{\partial U_{j}}=0$
and so $\frac{\partial g_{l}}{\partial V_{p}}(Z)\frac{\partial g_{p}(Z)}{\partial V_{q_{1}^{p}}}\cdots\frac{\partial g_{q_{m_{l}}^{p}}(Z)}{\partial U_{j}}=0$
in contradiction with the existence of $\left(q_{1}^{l},\ldots q_{m_{l}+1}^{l}\right)$,
thus $\frac{\partial  \tilde V_{l,i+m_{l}+2}}{\partial  \tilde U_{j,i}}\neq0.$ We
proved that if $\mathcal{P}(m)$ is true for sequence $m_{l}\leq m$,
the property holds as well for sequences of length $m_{l}+1$ and
so $\mathcal{P}(m+1)$ is true as well. This concludes the recursion
and the proof.
\end{proof}

\section*{Appendix 2: $\overline{u}_{Z_{0}}$ and $\widehat{Z_{0}}$ formal
derivation}

Let us introduce the following matrices for $i=0,\ldots,n$: 
\[
\mathbf{\overline{A}}_{\theta,i}=\left(\begin{array}{cc}
\mathbf{A}_{\theta,i} & \Delta_{i}r_{\theta}(t_{i})\\
0_{1,d} & 1
\end{array}\right),\,Q_{i}=\left(\begin{array}{cc}
C^{T}C & -C^{T}Y_{i}\\
-Y_{i}^{T}C & Y_{i}^{T}Y_{i}
\end{array}\right),\boldsymbol{\varGamma}_{\sigma}\left(t_{i}\right)=\left(\begin{array}{c}
\sqrt{\triangle }\varGamma_{\sigma}(t_{i})\\
0_{1,d_{u}}
\end{array}\right),P_{i}=\frac{1}{w}I_{m}.
\]
Let us extend the process $Z$ by introducing $\mathbf{Z}_{0}^{d}=\left(Z_{0}^{T},1\right)^{T}$ and by iteration, $\mathbf{Z}_{i}^{d}=\left(\tilde Z_{i},1\right)$ the extended state-space
variable:
\begin{equation}
\left\{ \begin{array}{l}
\mathbf{Z}_{i+1}^{d}=\mathbf{\overline{A}}_{\theta,i}\mathbf{Z}_{i+1}^{d}+\boldsymbol{\varGamma}_{\sigma}\left(t_{i}\right)u_{i}
\end{array}\right.\label{eq:extended_discretized_equation}
\end{equation}
as well as the reformulation of the control problem for  the extended process $\mathbf{Z}^d$:
\begin{equation}
\begin{array}{ll}
\textrm{Minimize in \ensuremath{u}:} & C_{w}(u|Y;Z_{0})=\left(\mathbf{Z}_{n}^{d}\right)^{T}Q_{n}\mathbf{Z}_{n}^{d}+\sum_{i=1}^{n-1}\left(\left(\mathbf{Z}_{i}^{d}\right)^{T}Q_{i}\mathbf{Z}_{i}^{d}+u_{i}^{T}P_{i}u_{i}\right)\\
\textrm{Subject to: } & \left\{ \begin{array}{l}
\mathbf{Z}_{i+1}^{d}=\mathbf{\overline{A}}_{\theta,i}\mathbf{Z}_{i+1}^{d}+\boldsymbol{\varGamma}_{\sigma}\left(t_{i}\right)u_{i}\\
\mathbf{Z}_{0}^{d}=\left(Z_{0},\,1\right)
\end{array}\right.
\end{array}\label{eq:extended_OC_Cost_Regulator_Form}
\end{equation}
 By introducing $E_{\theta,k}^{d}$ the solution of the discrete
Riccati equation: 
\begin{equation}
\left\{ \begin{array}{l}
R_{i}=Q_{i}+\mathbf{\overline{A}}_{\theta,i}^{T}R_{i+1}\mathbf{\overline{A}}_{\theta,i}\\
-\mathbf{\overline{A}}_{\theta,i}^{T}R_{i+1}\boldsymbol{\varGamma}_{\sigma}(t_{i})(P_{i}+\boldsymbol{\varGamma}_{\sigma}(t_{i}){}^{T}R_{i+1}\boldsymbol{\varGamma}_{\sigma}(t_{i}))^{-1}\boldsymbol{\varGamma}_{\sigma}(t_{i}){}^{T}R_{j+1}\mathbf{\overline{A}}_{\theta,i}\\
R_{n}=Q_{n}
\end{array}\right.\label{eq:extended_Riccati_equation}
\end{equation}
we know that the solution of (\ref{eq:extended_OC_Cost_Regulator_Form})
is equal to
\[
\overline{u}_{Z_{0},i}=-(P_{i}+\boldsymbol{\varGamma}_{\sigma}(t_{i}){}^{T}R_{i+1}\boldsymbol{\varGamma}_{\sigma}(t_{i}))^{-1}\boldsymbol{\varGamma}_{\sigma}(t_{i}){}^{T}R_{i+1}\mathbf{\overline{A}}_{\theta,i}\overline{\mathbf{Z}_{i}^{d}}.
\]
The corresponding cost value is given by $C_{w}(\overline{u}_{Z_{0}}|Y;Z_{0})=\left(\mathbf{Z}_{0}^{d}\right)^{T}R_{0}\overline{\mathbf{Z}_{0}^{d}}$
(see \cite{Bertsekas2005}). We now derive by reversed time induction
that the solution of (\ref{eq:extended_Riccati_equation}) is a sequence
of symmetric matrices i.e. $R_{i}$ has the block structure $R_{i}^{d}=\left(\begin{array}{cc}
E_{i}& h_{i}\\
\left(h_{i}\right)^{T} & \alpha_{i}^{d}
\end{array}\right)$ where each block respects a specific finite difference equation.
Let us assume the decomposition holds for $i+1$ i.e. $R_{i+1}=\left(\begin{array}{cc}
E_{i+1} & h_{i+1}\\
\left(h_{i+1}\right)^{T} & \alpha_{i+1}^{d}
\end{array}\right)$ 
Let us compute $\mathbf{\overline{A}}_{\theta,i}^{T}R_{i+1}\mathbf{\overline{A}}_{\theta,i}$. We have: 
\[
\begin{array}{l}
\mathbf{\overline{A}}_{\theta,i}^{T}R_{i+1}\mathbf{\overline{A}}_{\theta,i}\\
=\left(\begin{array}{cc}
\mathbf{A}_{\theta,i} & \Delta_{i}r_{\theta}(t_{i})\\
0_{1,d} & 1
\end{array}\right)^{T}\left(\begin{array}{cc}
E_{i+1} & h_{i+1}\\
\left(h_{i+1}\right)^{T} & \alpha_{i+1}^{d}
\end{array}\right)\left(\begin{array}{cc}
\mathbf{A}_{\theta,i} & \Delta_{i}r_{\theta}(t_{i})\\
0_{1,d} & 1
\end{array}\right)\\
=\left(\begin{array}{cc}
\mathbf{A}_{\theta,i} & \Delta_{i}r_{\theta}(t_{i})\\
0_{1,d} & 1
\end{array}\right)^{T}\left(\begin{array}{lc}
E_{i+1}\mathbf{A}_{\theta,i} & \Delta_{i}E_{i+1}r_{\theta}(t_{i})+h_{i+1}\\
\left(h_{i+1}\right)^{T}\mathbf{A}_{\theta,i} & \Delta_{i}\left(h_{i+1}\right)^{T}r_{\theta}(t_{i})+\alpha_{i+1}^{d}
\end{array}\right)\\
=\left(\begin{array}{cc}
\mathbf{A}_{\theta,i}^{T} & 0_{d,1}\\
\Delta_{i}r_{\theta}(t_{i})^{T} & 1
\end{array}\right)\left(\begin{array}{lc}
E_{i+1}\mathbf{A}_{\theta,i} & \Delta_{i}E_{i+1}r_{\theta}(t_{i})+h_{i+1}\\
\left(h_{i+1}\right)^{T}\mathbf{A}_{\theta,i} & \Delta_{i}\left(h_{i+1}\right)^{T}r_{\theta}(t_{i})+\alpha_{i+1}^{d}
\end{array}\right)\\
=\left(\begin{array}{lc}
\mathbf{A}_{\theta,i}^{T}E_{i+1}\mathbf{A}_{\theta,i} & \Delta_{i}\mathbf{A}_{\theta,i}^{T}E_{i+1}r_{\theta}(t_{i})+\mathbf{A}_{\theta,i}^{T}h_{i+1}\\
\Delta_{i}r_{\theta}(t_{i})^{T}E_{i+1}\mathbf{A}_{\theta,i}+\left(h_{i+1}\right)^{T}\mathbf{A}_{\theta,i} & \Delta_{i}r_{\theta}(t_{i})^{T}\left(\Delta_{i}E_{i+1}r_{\theta}(t_{i})+2h_{i+1}\right)+\alpha_{i+1}^{d}
\end{array}\right)\\
=F(E_{i+1}):=\left(\begin{array}{cc}
F_{1}(E_{i+1}) & F_{2}(E_{i+1},h_{i+1})\\
F_{2}(E_{i+1},h_{i+1})^{T} & F_{3}(E_{i+1},h_{i+1},\alpha_{\theta,k+1}^{d})
\end{array}\right)
\end{array}
\]
where the symmetric matrix $F$ is easily derivable by identification
and the $d$ dimensional square matrix $F_{1}$ only depends on $E_{i+1}$.
Now let us compute the expression in $[P_{i}+\boldsymbol{\varGamma}_{\sigma}(t_{i}){}^{T}R_{i+1}\boldsymbol{\varGamma}_{\sigma}(t_{i})]^{-1}$. We obtain:

\[
\begin{array}{lll}
[P_{i}+\boldsymbol{\varGamma}_{\sigma}(t_{i}){}^{T}R_{i+1}\boldsymbol{\varGamma}_{\sigma}(t_{i})]^{-1} & = & \left[\frac{1}{w}I_{m}+\left(\begin{array}{c}
\sqrt{\triangle }\varGamma_{\sigma}(t_{i})\\
0_{1,d_{u}}
\end{array}\right)^{T}\left(\begin{array}{cc}
E_{i+1} & h_{i+1}\\
\left(h_{i+1}\right)^{T} & \alpha_{i+1}^{d}
\end{array}\right)\left(\begin{array}{c}
\sqrt{\triangle }\varGamma_{\sigma}(t_{i})\\
0_{1,d_{u}}
\end{array}\right)\right]^{-1}\\
 & = & \left[\frac{1}{w}I_{m}+\triangle \left(\begin{array}{c}
E_{i+1}\varGamma_{\sigma}(t_{i})\\
\left(h_{i+1}\right)^{T}\varGamma_{\sigma}(t_{i})
\end{array}\right)\right]^{-1}\\
 & = & \left[\frac{1}{w}I_{m}+\triangle \left(\begin{array}{cc}
\varGamma_{\sigma}(t_{i})^{T} & 0_{d_{u},1}\end{array}\right)\left(\begin{array}{c}
E_{i+1}\varGamma_{\sigma}(t_{i})\\
\left(h_{i+1}\right)^{T}\varGamma_{\sigma}(t_{i})
\end{array}\right)\right]^{-1}\\
 & = & \left[\frac{1}{w}I_{m}+\triangle \varGamma_{\sigma}(t_{i})^{T}E_{i+1}\varGamma_{\sigma}(t_{i})\right]^{-1}:=G(E_{i+1})
\end{array}
\]
with $G$ a symmetric matrix valued function (as the inverse of a
symmetric square matrix). Moreover, we have:

\[
\begin{array}{l}
\begin{array}{l}
\mathbf{\overline{A}}_{\theta,i}^{T}R_{i+1}\boldsymbol{\varGamma}_{\sigma}(t_{i})\\
=\left(\begin{array}{cc}
\mathbf{A}_{\theta,i} & \Delta_{i}r_{\theta}(t_{i})\\
0_{1,d} & 1
\end{array}\right)^{T}\left(\begin{array}{cc}
E_{i+1} & h_{i+1}\\
\left(h_{i+1}\right)^{T} & \alpha_{i+1}^{d}
\end{array}\right)\left(\begin{array}{c}
\sqrt{\triangle }\varGamma_{\sigma}(t_{i})\\
0_{1,d_{u}}
\end{array}\right)\\
=\sqrt{\triangle }\left(\begin{array}{cc}
\mathbf{A}_{\theta,i}^{T} & 0_{d,1}\\
\Delta_{i}r_{\theta}(t_{i})^{T} & 1
\end{array}\right)\left(\begin{array}{c}
E_{i+1}\varGamma_{\sigma}(t_{i})\\
\left(h_{i+1}\right)^{T}\varGamma_{\sigma}(t_{i})
\end{array}\right)\\
=\sqrt{\triangle }\left(\begin{array}{c}
\mathbf{A}_{\theta,i}^{T}E_{i+1}\varGamma_{\sigma}(t_{i})\\
\Delta_{i}r_{\theta}(t_{i})^{T}E_{i+1}\varGamma_{\sigma}(t_{i})+\left(h_{i+1}\right)^{T}\varGamma_{\sigma}(t_{i})
\end{array}\right)
\end{array}\end{array}
\]
so we can compute

\[
\begin{array}{l}
\mathbf{\overline{A}}_{\theta,i}^{T}R_{i+1}\boldsymbol{\varGamma}_{\sigma}(t_{i})G(E_{i+1})\boldsymbol{\varGamma}_{\sigma}(t_{i}){}^{T}R_{j+1}\mathbf{\overline{A}}_{\theta,i}\\
=\sqrt{\triangle }\left(\begin{array}{c}
\mathbf{A}_{\theta,i}^{T}E_{i+1}\varGamma_{\sigma}(t_{i})\\
\Delta_{i}r_{\theta}(t_{i})^{T}E_{i+1}\varGamma_{\sigma}(t_{i})+\left(h_{i+1}\right)^{T}\varGamma_{\sigma}(t_{i})
\end{array}\right)G(E_{i+1})\sqrt{\triangle }\left(\begin{array}{c}
\mathbf{A}_{\theta,i}^{T}E_{i+1}\varGamma_{\sigma}(t_{i})\\
\Delta_{i}r_{\theta}(t_{i})^{T}E_{i+1}\varGamma_{\sigma}(t_{i})+\left(h_{i+1}\right)^{T}\varGamma_{\sigma}(t_{i})
\end{array}\right)^{T}\\
=\triangle \left(\begin{array}{c}
\mathbf{A}_{\theta,i}^{T}E_{i+1}\varGamma_{\sigma}(t_{i})\\
\left(\Delta_{i}E_{i+1}r_{\theta}(t_{i})+h_{i+1}\right)^{T}\varGamma_{\sigma}(t_{i})
\end{array}\right)G(E_{i+1})\left(\begin{array}{cc}
\varGamma_{\sigma}(t_{i})^{T}E_{i+1}\mathbf{A}_{\theta,i} & \varGamma_{\sigma}(t_{i})^{T}\left(h_{i+1}+\Delta_{i}E_{i+1}r_{\theta}(t_{i})\right)\end{array}\right)\\
:=\triangle \left(\begin{array}{cc}
H_{1}(E_{i+1}) & H_{2}(E_{i+1},h_{i+1})\\
H_{2}(E_{i+1},h_{i+1})^{T} & H_{3}(E_{i+1},h_{i+1})
\end{array}\right).
\end{array}
\]
By re-injecting all the derived expression in (\ref{eq:extended_Riccati_equation}),
we obtain:

\[
\begin{array}{l}
R_{i}=\left(\begin{array}{cc}
F_{1}(E_{i+1}) & F_{2}(E_{i+1},h_{i+1})\\
F_{2}(E_{i+1},h_{i+1})^{T} & F_{3}(E_{i+1},h_{i+1},\alpha_{\theta,k+1}^{d})
\end{array}\right)+\left(\begin{array}{cc}
C^{T}C & -C^{T}Y_{i}\\
-Y_{i}^{T}C & Y_{i}^{T}Y_{i}
\end{array}\right)\\
-\triangle \left(\begin{array}{cc}
H_{1}(E_{i+1}) & H_{2}(E_{i+1},h_{i+1})\\
H_{2}(E_{i+1},h_{i+1})^{T} & H_{3}(E_{i+1},h_{i+1})
\end{array}\right)
\end{array}
\]
Thus $R_{i}$ is indeed symmetric and has the required
form. Hence the recursion. We also obtain the following finite difference
equation

\[
\begin{array}{l}
E_{i}^{d}=\mathbf{A}_{\theta,i}^{T}E_{i+1}\mathbf{A}_{\theta,i}+C^{T}C-\triangle \mathbf{A}_{\theta,i}^{T}E_{i+1}\varGamma_{\sigma}(t_{i})G(E_{i+1})\varGamma_{\sigma}(t_{i})^{T}E_{i+1}\mathbf{A}_{\theta,i}\\
h_{i}=\Delta_{i}\mathbf{A}_{\theta,i}^{T}E_{i+1}r_{\theta}(t_{i})+\mathbf{A}_{\theta,i}^{T}h_{i+1}-C^{T}Y_{i}\\
-\triangle \mathbf{A}_{\theta,i}^{T}E_{i+1}\varGamma_{\sigma}(t_{i})G(E_{i+1})\varGamma_{\sigma}(t_{i})^{T}\left(h_{i+1}+\Delta_{i}E_{i+1}r_{\theta}(t_{i})\right).
\end{array}
\]
From the previous expression for $R_{i}$, we derive for the related
optimal control $\overline{u_{\theta}^{d}}$:
\[
\begin{array}{l}
\overline{u}_{i}=-(P_{i}+\boldsymbol{\varGamma}_{\sigma}(t_{i}){}^{T}R_{i+1}\boldsymbol{\varGamma}_{\sigma}(t_{i}))^{-1}\boldsymbol{\varGamma}_{\sigma}(t_{i}){}^{T}R_{i+1}\mathbf{\overline{A}}_{\theta,i}\left(\begin{array}{l}
\overline{Z}_{i}\\
1
\end{array}\right)\\
=-G(E_{i+1})\sqrt{\triangle }\left(\begin{array}{cc}
\varGamma_{\sigma}(t_{i})^{T} & 0_{d_{u},1}\end{array}\right)\left(\begin{array}{cc}
E_{i+1} & h_{i+1}\\
\left(h_{i+1}\right)^{T} & \alpha_{i+1}^{d}
\end{array}\right)\left(\begin{array}{cc}
\mathbf{A}_{\theta,i} & \Delta_{i}r_{\theta}(t_{i})\\
0_{1,d} & 1
\end{array}\right)\left(\begin{array}{l}
\overline{Z}_{i}\\
1
\end{array}\right)\\
=-\sqrt{\triangle }G(R_{\theta,k+1}^{d})\left(\begin{array}{cc}
\varGamma_{\sigma}(t_{i})^{T} & 0_{d_{u},1}\end{array}\right)\left(\begin{array}{cc}
E_{i+1} & h_{i+1}\\
\left(h_{i+1}\right)^{T} & \alpha_{i+1}^{d}
\end{array}\right)\left(\begin{array}{c}
\mathbf{A}_{\theta,i}\overline{Z}_{i}+\Delta_{i}r_{\theta}(t_{i})\\
1
\end{array}\right)\\
=-\sqrt{\triangle }G(R_{\theta,k+1}^{d})\left(\begin{array}{cc}
\varGamma_{\sigma}(t_{i})^{T} & 0_{d_{u},1}\end{array}\right)\left(\begin{array}{c}
E_{i+1}\left(\mathbf{A}_{\theta,i}\overline{Z}_{i}+\Delta_{i}r_{\theta}(t_{i})\right)+h_{i+1}\\
\left(h_{i+1}\right)^{T}\left(\mathbf{A}_{\theta,i}\overline{Z}_{i}+\Delta_{i}r_{\theta}(t_{i})\right)+\alpha_{i+1}^{d}
\end{array}\right)\\
=-\sqrt{\triangle }G(R_{\theta,k+1}^{d})\varGamma_{\sigma}(t_{i})^{T}\left(E_{i+1}\left(\mathbf{A}_{\theta,i}\overline{Z}_{i}+\Delta_{i}r_{\theta}(t_{i})\right)+h_{i+1}\right).
\end{array}
\]
From these equations, we easily derive that
\[
\begin{array}{l}
\widehat{Z_{0}}=\arg\min_{Z_{0}}C_{w}(\overline{u}_{Z_{0}}|Y;Z_{0})\\
=\arg\min_{Z_{0}}\left\{ \left(\mathbf{Z}_{0}^{d}\right)^{T}\left(\begin{array}{cc}
E_{0}^{d} & h_{0}^{d}\\
\left(h_{0}^{d}\right)^{T} & \alpha_{0}^{d}
\end{array}\right)\overline{\mathbf{Z}_{0}^{d}}\right\} \\
=\arg\min_{Z_{0}}\left\{ Z_{0}^{T}R_{0}^{d}Z_{0}+2\left(h_{0}^{d}\right)^{T}Z_{0}+\alpha_{0}^{d}\right\} \\
=-\left(E_{0}^{d}\right)^{-1}h_{0}^{d}.
\end{array}
\]

\bibliographystyle{Chicago}

\bibliography{biblio_hypo_sde}